\input ./style/arxiv-general.cfg
\documentclass[aos,MSNbibl,seceqn,dvips]{arximspdf}
\makeatletter
   \@ifpackageloaded{graphicx}{}{\usepackage{graphicx}}
\makeatother

\doi{10.1214/15-AOS1334}
\volume{43}
\issue{5}
\pubyear{2015}
\firstpage{1986}
\lastpage{2018}
\docsubty{FLA}

\makeatletter
\newcommand{\eqref}[1]{(\ref{#1})}
\newcommand{\mc}{\mathop{\mathrm{ mc}} }
\newcommand{\argmin}{\mathop{\arg\min}}

\newtheorem{thmm}{Theorem}
\newtheorem{lem}{Lemma}
\newproclaim{rem}{Remark}

\newtheorem{cor}{Corollary}
\newproclaim{assump}{Assumption}
\newproclaim{definition}{Definition}[section]
\newproclaim{exam}{Example}

\newcommand{\ola}{\overline{\lambda}}
\newcommand{\E}{\mathbb{E}}
\renewcommand{\P}{\mathbb{P}}

\newcommand{\Tr}{\operatorname{tr}}
\newcommand{\TV}{\operatorname{TV}}
\newcommand{\Beta}{\operatorname{B}}
\newcommand{\pen}{\operatorname{pen}}

\newcommand{\leqa}{\lesssim}
\newcommand{\geqa}{\gtrsim}

\newcommand{\ra}{\rightarrow}{}
\newcommand{\rn}{\sqrt{n}}

\newcommand{\al}{\alpha}
\newcommand{\be}{\beta}
\newcommand{\ep}{\varepsilon}
\newcommand{\Ga}{\Gamma}
\newcommand{\la}{\lambda}
\newcommand{\si}{\sigma}

\newcommand{\cN}{\mathcal{N}}
\newcommand{\cS}{\mathcal{S}}
\newcommand{\cT}{\mathcal{T}}
\newcommand{\cV}{\mathcal{V}}
\newcommand{\cY}{\mathcal{Y}}

\newcommand{\RR}{\mathbb{R}}

\newcommand{\Given}{ | }

\makeatother

\begin{document}
\begin{frontmatter}

\title{Bayesian linear regression with sparse priors}
\runtitle{Bayes sparse linear regression}

\begin{aug}
\author[A]{\fnms{Isma\"el}~\snm{Castillo}\thanksref{T1,m1}\ead[label=e1]{ismael.castillo@math.cnrs.fr}},
\author[B]{\fnms{Johannes}~\snm{Schmidt-Hieber}\thanksref{T2,m2}\ead[label=e2]{schmidthieberaj@math.leidenuniv.nl}}\\
\and
\author[B]{\fnms{Aad}~\snm{van der Vaart}\corref{}\thanksref{T2,m2}\ead[label=e3]{avdvaart@math.leidenuniv.nl}}
\runauthor{I. Castillo, J. Schmidt-Hieber and A. van der Vaart}
\thankstext{T1}{Supported in part by ANR Grants ``Banhdits''
ANR-2010-BLAN-0113-03 and ``Calibration'' ANR-2011-BS01-010-01.}
\thankstext{T2}{Supported by the European Research Council under ERC
Grant Agreement 320637.}

\affiliation{CNRS Paris\thanksmark{m1} and Leiden
University\thanksmark{m2}}

\address[A]{I. Castillo\\
CNRS - LPMA\\
Universit\'es Paris VI \& VII \\
B\^atiment Sophie Germain\\
75205 Paris Cedex 13\\
France\\
\printead{e1}}

\address[B]{J. Schmidt-Hieber\\
A. van der Vaart\\
Mathematical Institute\\
Leiden University\\
Niels Bohrweg 1 \\
2333 CA Leiden\\
The Netherlands\\
\printead{e2}\\
\phantom{E-mail:\ }\printead*{e3}}
\end{aug}

%
\received{\smonth{3} \syear{2014}}
%
\revised{\smonth{3} \syear{2015}}

%
\begin{abstract}
We study full Bayesian procedures for high-dimensional linear
regression under sparsity
constraints. The prior is a mixture of point masses at zero and
continuous distributions. Under
compatibility conditions on the design matrix, the posterior
distribution is shown to
contract at the optimal rate for recovery of the unknown sparse vector,
and to give optimal
prediction of the response vector. It is also shown to select the
correct sparse model, or at
least the coefficients that are significantly different from zero. The
asymptotic shape of the
posterior distribution is characterized and employed to the
construction and study of credible
sets for uncertainty quantification.
\end{abstract}

%

%
\end{frontmatter}

\section{Introduction} \label{sec.intro}

Consider estimation of a parameter $\be\in\RR^p$ in the linear
regression model
%
\begin{equation}
\label{Model} Y= X\be+\ep,
\end{equation}
where $X$ is a given, deterministic $(n\times p)$ matrix, and $\ep$ is
an $n$-variate standard
normal vector. The model is standard, but we are interested in the
\emph{sparse} setup, where $n\leq
p$, and possibly $n\ll p$, and ``many'' or ``most'' of the coefficients
$\be_i$ of the parameter vector
are zero, or close to zero. We study a Bayesian approach based on
priors that set a selection of
coefficients $\be_i$ a priori to zero; equivalently, priors that
distribute their mass over models
that use only a (small) selection of the columns of $X$. Bayes's
formula gives a posterior
distribution as usual. We study this under the ``frequentist''
assumption that the data $Y$ has in
reality been generated according to a given (sparse) parameter $\be
^0$. The expectation under the previous distribution is denoted
$\mathbb{E}_{\be^0}$.

Specifically, we consider a prior $\Pi$ on $\be$ that first selects a
\emph{dimension}
$s$ from a prior $\pi_p$ on the set $\{0,\ldots,p\}$,
next a random subset $S\subset\{1,2,\ldots, p\}$ of cardinality
$|S|=s$ and finally a set of
nonzero values $\be_S:=\{\be_i: i\in S\}$ from a prior density $g_S$
on $\RR^S$.
Formally, the prior on $(S,\be)$ can be expressed as
%
\begin{equation}
\label{DefPrior} (S,\be)\mapsto{\pi_p\bigl(|S|\bigr)}\frac{1}{{p\choose|S|}}
g_S(\be_S) \delta _0(\be_{S^c}),
\end{equation}
where the term $\delta_0(\be_{S^c})$
refers to the coordinates $\be_{S^c}:=(\be_i: i\in S^c)$
being zero. We focus on the situation where $g_S$ is a product $\otimes
g$ of
densities over the coordinates in~$S$, for $g$ a fixed continuous
density on $\RR$,
with the Laplace density as an important special case.
The prior $\pi_p$ is crucial for expressing the ``sparsity'' of the parameter.
One of the main findings of this paper is that weights $\pi_p(s)$ that
decrease slightly faster than exponential in the dimension $s$ give
good performance.

Priors of the type of \eqref{DefPrior} were considered by many authors,
including \cite{mitchellbeauchamp,georgefoster,george,ishwaranrao,yuanglin,ScottBerger10,BottoloRichardson}. Other related
contributions include \cite{abragrin10,ariasclounici,2014arXiv1406.7718M}. The paper \cite{spa} contains
a theoretical analysis similar to the present paper, but restricted to the
special case that the regression
matrix $X$ is the identity and $p=n$; see Example~\ref{ExampleSequenceModel}.
The general model \eqref{Model} shares some features with this special
case, but is different in
that it must take account of the noninvertibility
of $X$ and its interplay with the sparsity assumption, especially for
the case of recovering the
parameter $\be$, as opposed to estimating the mean $X\be$. While the
proofs in \cite{spa}
use a factorization of the model along the coordinate axes, exponential tests
and entropy bounds, in the present paper we employ a direct and refined
analysis of the
posterior ratio \eqref{bayes}, exploiting the specific form of the
prior Laplace density $g$.
Furthermore, even for the case that $X$ is the identity matrix, the
present paper provides
several new results of interest: distributional
approximations to the posterior distribution, insight in the
scaling of the prior on the nonzero coordinates and oracle formulations
of the contraction rates.

Algorithms for the computation of the posterior distribution corresponding
to~\eqref{DefPrior}, especially for the ``spike and slab'' prior
described in Example~\ref{example.slabspike} below,
are routine for small dimensions $p$ and $n$
(e.g., \cite{mitchellbeauchamp,georgefoster,ishwaranrao,Dellaportas,yuanglin,ScottBerger10,JiSchnmidler,HansDobraWest}).
For large dimensions the resulting computations are intensive, due to
the large number of possible
submodels $S$. Many authors are currently developing algorithms that
can cope with larger numbers
of covariates, in the sparse setup considered in the present paper. In
Section~\ref{SectionComputationalAspects} we review recent progress
on various methods, of which some are feasible for values of $p$ up to
hundreds or thousands \cite
{StingoVannucci,LiZhang,BottoloRichardson,BottoloRichardson2,RichardsonBottolo,ShiDunson,SchaferChopin,2013arXiv1312.5658S,RockovaGeorge,Ormerod,spa}.
Although this upper bound will increase in the coming years, clearly it
falls far short of
the dimensions attainable by (point) estimation methods based on convex
programming,
such as the LASSO. Other Bayesian approaches to sparse regression that
do not explicitly include
model selection (e.g., \cite{griffinbrown,carvalhopolsonscott,hans})
can cope with somewhat higher dimensions,
but truly high-dimensional models are out of reach of fully Bayesian
methods at the present time.

Not surprisingly to overcome the nonidentifiability of the full
parameter vector $\be$
in the overspecified model \eqref{Model}, we borrow from the work on
sparse regression within the non-Bayesian framework; see \cite{don06,fdr06,candestao07,BRT09,buehlmannvdGeer2011,rigollet-tsybakov11,KoltchinskiiLouniciTsybakov,zhang10,zhanghuang,buhlmannvandegeerzhang}. Good performance of
the posterior distribution is shown under \emph{compatibility} and
\emph{smallest sparse
eigenvalue} conditions; see Section~\ref{sec.recovery}. Although the
constants in these
results are not as sharp as results for the LASSO, the posterior
contraction rates obtained
are broadly comparable to convergence rates of the LASSO.

The LASSO and its variants are important frequentist methods for sparse
signal recovery. As the
LASSO is a posterior mode (for an i.i.d. Laplace prior on the $\be
_i$), it may seem to give an
immediate link between Bayesian and non-Bayesian methods. However, we
show in
Section~\ref{sec.spars_and_LASSO} that the LASSO is essentially
non-Bayesian, in the sense that the
corresponding \emph{full} posterior distribution is a useless object.

In contrast, the posterior distribution
resulting from the prior \eqref{DefPrior} gives both reasonable
reconstruction of the parameter $\be$ and
a quantification of uncertainty through the spread in the posterior
distribution.
We infer this from combining results on the contraction rate of the
full posterior distribution with distributional approximations. The
latter show that
the posterior distribution behaves asymptotically as
a mixture of Bernstein--von Mises type approximations to submodels,
where the
location of the mixture components depends on the setting.
The latter approximations are new, also for the special case that $X$
is the identity matrix.

It is crucial for these results that the prior \eqref{DefPrior} models
sparsity through the
\emph{model selection} prior $\pi_p$, and separates this from
modeling the nonzero coordinates
through the prior densities $g_S$. For instance, in the case that $g_S$
is a product of Laplace
densities, this allows the scale parameter to be constant or even to
tend to zero, thus making this
prior uninformative. This is in stark contrast to the choice of the
smoothing parameter in the
(Bayesian) LASSO, which must tend to infinity in order to shrink
parameters to zero, where it cannot
differentiate between truly small and nonzero parameters. Technically
this has the consequence that
the essential part of the proofs is to show that the posterior
distribution concentrates on sets of
small dimension. This sets it apart from the frequentist literature on
sparse regression,
although, as mentioned, many essential ideas reappear here in a
Bayesian framework.

The paper is organized as follows. In Section~\ref{sec.recovery} we
present the main results of the
paper. We specialize to Laplace priors on the nonzero coefficients and
investigate the ability of
the posterior distribution to recover the parameter vector $\be$, the
predictive vector $X\be$
and the set of nonzero coordinates. Furthermore, we derive a
distributional approximation to the
posterior distribution, and apply this to construct and study credible
sets. In
Section~\ref{sec.spars_and_LASSO} we present the negative result on
the Bayesian interpretation of
the LASSO. Next in Section~\ref{SectionArbitraryDesign} we show that
for recovery
of only the predictive vector $X\be$, significantly milder conditions
than in
Section~\ref{sec.recovery} suffice. Proofs are deferred to
Section~\ref{SectionProofs}
and the supplementary material \cite{Supplement}.

\subsection{Notation}
For a vector $\be\in\RR^p$ and a set $S\subset\{1,2,\ldots,p\}$ of indices,
$\be_S$ is the vector $(\be_i)_{i\in S}\in\RR^S$, and $|S|$ is the
cardinality of $S$.
The \emph{support} of the parameter $\be$ is the set $S_\be=\{i: \be
_i\neq0\}$.
The support of the true parameter $\be^0$ is denoted $S_0$,
with cardinality $s_0:=|S_0|$. Similarly, for a generic vector $\be
^*$, we
write $S_*=S_{\be^*}$ and $s_*=|S_*|$. We write $s=|S|$ if there is no
ambiguity to
which set $S$ is referred to.
For $1\le q< \infty$ and $\be\in\RR^p$, let $\|\be\|_q:=(\sum_{i=1}^p |\be_i|^q)^{1/q}$.


We let $X_{\cdot,i}$ be the $i$th column of $X$, and
%
\begin{equation}
\label{DefNormX} \|X\| 
=\max_{i=1,\ldots,p}\|X_{\cdot,i}
\|_2 =\max_{i=1,\ldots,p} \bigl(X^tX
\bigr)_{i,i}^{1/2}.
\end{equation}
For the prior $\Pi$ defined above, Bayes's formula gives the following
expression for the posterior distribution $\Pi[\cdot| Y]$. For any
Borel set $B$ of $\mathbb{R}^p$,
%
\begin{equation}
\label{bayes} \Pi[B| Y] = \int_B e^{-\|Y-X\be\|_2^2/2}\,d\Pi(
\be) \Big/ \int e^{-\|Y-X\be\|_2^2/2}\,d\Pi(\be).
\end{equation}


\section{Main results}
\label{sec.recovery}
\subsection{Prior}
In this section we consider the prior \eqref{DefPrior}, with $g_S$ the
product of $|S|$
Laplace densities $\be\mapsto2^{-1}\la\exp(-\la|\be|)$. We allow
the (inverse) scale parameter $\la$ to change with $p$, within the range,
with $\|X\|$ defined in \eqref{DefNormX},
%
\begin{equation}
\frac{\|X\|}{p} \leq\la\le2\ola,\qquad \ola=2\|X\|\sqrt{\log p}. \label{eq.lambda_cond}
\end{equation}
The quantity $\overline\la$ in the upper bound is the usual value of
the regularization parameter
$\la$ of the LASSO [as in \eqref{EqLASSO} below].
Its large value causes the LASSO to shrink many coordinates $\be_i$ to
zero, as is desired
in the sparse situation. However, in our Bayesian setup, sparsity
should be induced by
model selection, through the prior $\pi_p$ on the model dimension, and
the Laplace
prior densities model only the nonzero coordinates. Large values of
$\la$
would shrink the nonzero coordinates to zero, which is clearly
undesirable and unnatural.
Thus it is natural to assume $\la\ll\overline\la$, and fixed values
of $\la$,
and even values decreasing to zero, may well be natural, depending on
the regression setting.
We shall see that small values of $\la$ permit a distributional
approximation to
the posterior distribution centered at unbiased estimators. The results
below hold
for all $\lambda$ in the range \eqref{eq.lambda_cond}, but they
are meant to be
read for a specific sequence of $\lambda$ and are not suitable for
optimization over~$\la$.

The precise interpretation of the size of $\la$ is confounded with the
regression setting,
the error variance (which we have set to unity for simplicity of
notation) and the
scaling of the regression matrix. The following three special cases
shed some light on this.

\begin{exam}[(Sequence model)]
\label{ExampleSequenceModel}
In the \emph{sequence model} considered in \cite{JohnSilv04} and
\cite{spa}, the observation is a vector $(Y_1,\ldots, Y_n)$ of independent
coordinates $Y_i\sim N(\be_i,1)$. This corresponds to $X=I$ and $n=p$
in the present setting
\eqref{Model}, whence $\|X\|=1$. Condition \eqref{eq.lambda_cond}
then reduces to
$p^{-1}\le\la\le4\sqrt{\log p}$. Fixed values of $\la$, as
considered in
\cite{spa}, are easily included. As there is only one observation per
parameter,
it may not be unreasonable to consider $\la\ra0$, in order to create
noninformative priors for the nonzero coefficients. This is allowed
easily also.
\end{exam}

\begin{exam}[(Sequence model, multiple observations)]
\label{ExampleSequenceModelMultiple}
In an extension of the sequence model of the preceding example,
the $n$ observations are from normal distributions $N(\be_i,\sigma
_n^2)$ with variances
$\sigma_n^2\ra0$. By defining the $Y_i$ as $\sigma_n^{-1}$ times the original
observations, we can fit this into model \eqref{Model}, which has unit
error variances.
If we keep the original definition of the $\be_i$, then the regression
matrix is
$X=\sigma_n^{-1} I$, and hence
$\|X\|=\sigma_n^{-1}$. Condition \eqref{eq.lambda_cond} then reduces to
$\sigma^{-1}_n/ n\le\la\le4\sigma_n^{-1}\sqrt{\log n}$.
Fixed values of $\la$ are included if $n\sigma_n\gtrsim1$, and values
tending to zero if $n\sigma_n\ra\infty$. By sufficiency of the
sample mean
in the normal location problem this corresponds
to a sufficient number of replicate measurements on every parameter
$\be_i$
in the original problem.
\end{exam}

\begin{exam}[(Response model)]
\label{ExampleRegressionModel}
If every row of the regression equation $Y=X\be+\ep$ refers to a measurement
of an instance of a fixed relationship between an input vector
$X_{i,\cdot}\in\RR^p$
and the corresponding output $Y_i$, then the entry $X_{i,j}$ of $X$ is the
value of individual $i$ on the $j$th covariable. It is then reasonable
to think of
these entries as being sampled from some fixed distribution,
independent of $n$ and~$p$,
in which case $\|X\|$ will (typically) be of the order $\sqrt n$. A
fundamental example is the case where the entries of $X$ are
independent standard Gaussian $\cN(0,1)$.
Condition~\eqref{eq.lambda_cond} then reduces to
$\sqrt n/p\le\la\le4\sqrt n\sqrt{\log p}$. Fixed values of $\la$,
as considered in
\cite{spa}, are included, provided $p\gtrsim\sqrt n$.
\end{exam}

Although condition \eqref{eq.lambda_cond} does not exclude shrinkage
through large values of $\la$, as for the LASSO, the most interesting
situation is that sparsity is induced through model selection.
The prior $\pi_p$ on model dimension is crucial; it must downweight
big models,
but at the same time give sufficient mass to the true model.
Exponential decrease turns out to work.

\begin{assump}[(Prior dimension)]
There are constants $A_1, A_2, A_3, A_4>0$ with
%
\begin{equation}
A_1p^{-A_3}\pi_p(s-1)\le\pi_p(s)
\leq A_2p^{-A_4} \pi_p(s-1),\qquad s=1,\ldots,p.
\label{assump.on_the_dim_prior} 
\end{equation}
\end{assump}

\begin{exam}
Assumption \eqref{assump.on_the_dim_prior} is met by the priors of the form,
for constants $a, c>0$,
%
\begin{equation}
\pi_p(s) \propto c^{-s}p^{-as},\qquad s= 0,1,\ldots,p.
\label{eq.classofdimpriors}
\end{equation}
We refer to these priors as \emph{complexity priors}, as their rate of
decrease reflects the number of models ${p\choose s}$ of given size $s$
for $s\ll p$; cf. \cite{spa}.
\end{exam}

\begin{exam}[(Slab and spike)]\label{example.slabspike}
Modeling the coordinates $\be_1,\ldots,\be_p$ as i.i.d. variables
from a mixture $(1-r)\delta_0+rG$, of a Dirac measure $\delta_0$ at
zero and
a Laplace distribution $G$, is included in \eqref{DefPrior} with $\pi
_p$ the
binomial distribution with parameter $p$ and $r$. The size $r$ of the
point mass
at zero controls the model selection. The overall prior obtained
by choosing $r$ from a Beta $(1,p^u)$ hyper prior with $u>1$ satisfies
\eqref{eq.classofdimpriors}; cf. Example~2.2 in \cite{spa}. This
prior is universal in that it is free of unknown smoothing
parameters.
\end{exam}

To conclude the discussion on the prior, we briefly comment on the case
that the noise
vector has unknown variance $\sigma^2>0$; that is, we observe $Y=
X\beta+\sigma\varepsilon$. In
this case one may use an empirical Bayesian approach, which replaces
the unknown parameter
by an estimator, or a hierarchical Bayesian approach, which puts a
prior on $\sigma^2$, a common
choice being an inverse Gamma prior. Since $Y/\sigma= X(\beta/\sigma
) +\varepsilon$, it is natural to apply the prior, as in this paper,
to the parameter $\beta/\sigma$. Thus given $\sigma^2$ and a model
$S$, we choose the prior
density on the nonzero values $\beta_S=\{\beta_i:i\in S\}$ as the
product of $|S|$ Laplace
densities $\beta\mapsto\lambda/(2\sigma) \exp(-\lambda|\beta
|/\sigma)$, conditional
on the estimated or prior value of $\sigma$.

\subsection{Design matrix}
The parameter $\be$ in model \eqref{Model} is not estimable without
conditions on the regression matrix. For the interesting case $p>n$, it
is even necessarily
unidentifiable. If $\be$ is known to be sparse, then ``local
invertibility'' of the Gram
matrix $X^tX$ is sufficient for estimability, even in the case $p>n$.
This is made precise in the following definitions, which are based on
the literature, but with simplified notation suited to our Bayesian setup.
For accessibility we include short discussions on the relations between
the various concepts.

\begin{definition}[(Compatibility)]
\label{assump.comp}
The compatibility number of model $S\subset\{1,\ldots,p\}$ is given by
\[
\phi(S):=\inf \biggl\{\frac{\|X\be\|_2 |S|^{1/2}}{\|X\| \|\be_S\|_1 } : \|\be_{S^c}\|_1
\leq7 \|\be_S\|_1, \be_S\neq0 \biggr\}.
\]
\end{definition}

The compatibility number (which is $\phi_{\mathrm{comp}}(7,S)/\|X\|$ in the
notation of
\cite{buehlmannvdGeer2011}, page 157)
compares the $\ell^2$-norm of the predictive vector $X\be$
to the $\ell^1$-norm of the parameter $\be_S$. A model $S$ is
considered ``compatible''
if $\phi(S)>0$. It then satisfies the nontrivial inequality
$\|X\be\|_2|S|^{1/2}\ge\phi(S) \|X\| \|\be_S\|_1$. We shall see
that true vectors
$\be^0$ with compatible support $S_{\be^0}$ can be recovered from the data,
uniformly in a lower bound on the size of their compatibility numbers.

The number 7 has no particular interest, but for simplicity we use a
numerical value instead
of an unspecified constant. Since
the vectors $\be$ in the infimum satisfy $\|\be_S\|_1\le\|\be\|_1
\le8 \|\be_S\|_1$, it would not be a great loss of generality to
replace $\be_S$ in the denominator
of the quotient by $\be$. However, the factor $|S|^{1/2}$ in the numerator
may be seen as resulting from the comparison of the $\ell^1$- and
$\ell^2$-norms
of $\be_S$ through the Cauchy--Schwarz inequality: $\|\be_S\|_1\le
|S|^{1/2}\|\be_S\|_2$.
Replacing $\|\be_S\|_1/|S|^{1/2}$ by $\|\be_S\|_2$ would make the
compatibility number
smaller, and hence give a more restrictive condition.

The compatibility number involves the full vectors $\be$ (also their
coordinates
outside of $S$) and allows to reduce the recovery problem to sparse vectors.
The next two definitions concern sparse vectors only, but unlike the
compatibility
number, they are uniform in vectors up to a given dimension. In the
notation of
\cite{buehlmannvdGeer2011} (pages 156--157) the numbers in the
definitions are
the minima over $|S|\le s$ of the numbers $\Lambda_{\mathrm{min},1}(\Sigma
_{1,1}(S))/\|X\|$ and
$\Lambda_{\mathrm{min}}(\Sigma_{1,1}(S))/\|X\|$, respectively.

\begin{definition}[(Uniform compatibility in sparse vectors)]
\label{assump.uniform_comp}
The compatibility number in vectors of dimension $s$ is defined as
\[
\overline\phi(s):=\inf \biggl\{\frac{\|X\be\|_2 |S_\be|^{1/2}}{ \|
X\| \|\be\|_1 }: 0\neq|S_\be| \leq
s \biggr\}.
\]
\end{definition}

\begin{definition}[(Smallest scaled sparse singular value)]
\label{assump.smallest_sp_ev}
The smallest scaled singular value of dimension $s$ is defined as
%
\begin{equation}
\widetilde\phi(s) :=\inf \biggl\{\frac{\|X\be\|_2}{\|X\| \|\be\|
_2}: 0\neq|S_\be|
\leq s \biggr\}. \label{eq.def_widetilde_phi_K}
\end{equation}
\end{definition}

For recovery we shall impose that these numbers for $s$ equal to (a
multiple of) the dimension
of the true parameter vector are bounded away from zero.
Since $\|\be\|_1\le|S_\be|^{1/2}\|\be\|_2$ by the Cauchy--Schwarz
inequality,
it follows that $\widetilde\phi(s)\le\overline\phi(s)$, for any $s>0$.
The stronger assumptions on the design matrix imposed through
$\widetilde\phi(s)$ will
be used for recovery with respect to the $\ell^2$-norm, whereas the numbers
$\overline\phi(s)$ suffice for $\ell^1$-reconstruction. In
Definition~\ref{assump.smallest_sp_ev}, ``scaled'' refers to the
scaling of the matrix $X$ by division
by the maximum column length $\|X\|$; if the latter is unity, then
$\widetilde\phi(s)$ is just
the smallest scaled singular value of a submatrix of $X$ of dimension $s$.

The final and strongest invertibility condition is in terms of ``mutual
coherence'' of the regression
matrix, which is the maximum correlation between its columns.

\begin{definition}[(Mutual coherence)]
The \emph{mutual coherence number} is
\[
\mc(X)= \max_{1\le i\neq j\le p} \frac{|\langle X_{\cdot,i},
X_{\cdot,j}\rangle|}{\|X_{\cdot,i}\|_2 \|X_{\cdot,j}\|_2}.
\]
\end{definition}

The ``\emph{$(K,s)$ mutual coherence condition}'' is that this number
is bounded above by $(Ks)^{-1}$,
in which case reconstruction is typically possible for true vectors
$\be$ of dimension up to $s$.
As correlations are easy to interpret, conditions of this type, which
go back to \cite{don06},
have been used by many authors. (Notably, Bunea, Tsybakov and Wegkamp
\cite{BuneaTsybakovWegkamp} show that
for reconstructions using the $\ell^1$- and $\ell^2$-norms, taking
the maximum
over all correlations can be relaxed to a maximum over
pairs that involve at least one ``active'' coordinate.)
The following lemma shows that they are typically stronger than
conditions in terms
of compatibility numbers or sparse singular values.
The lemma is embodied in Lemma~2 in \cite{lou},
and is closely related to the inequalities obtained in \cite{vandegeer2009}.
For ease of reference we provide a proof in the supplementary material \cite{Supplement}.

\begin{lem}
\label{LemmaCoherence}
$\phi(S)^2 \ge \overline\phi(1)^2 - 15|S| \mc(X)$; $\overline\phi(s)^2 \ge \widetilde\phi(s)^2 \ge \overline\phi
(1)^2 - s \mc(X)$.
\end{lem}

By evaluating the infimum in
Definition~\ref{assump.uniform_comp} with $\be$ equal to unit
vectors, we
see that $\widetilde\phi(1)=\overline\phi(1)=\min_i\|X_{\cdot,i}\|_2/\|
X\|$, which will typically be
bounded away from zero. Thus the lemma implies that compatibility
numbers and sparse singular values
are certainly bounded
away from zero for models up to size a multiple of $1/\mc(X)$, that
is, models of size satisfying
the ``mutual coherence condition.'' This makes the mutual coherence the
strongest of the
three ``sparse invertibility'' indices introduced previously. We note
that the reverse
inequalities do not hold in general, and indeed the compatibility
constant can
easily be bounded away from zero, even if the mutual coherence number
is much larger than $1/s$.

For many other possible indices (including ``restricted isometry'' and
``irrepresentability''),
and extensive discussion of their relationships, we refer to
Sections~6.13 and~7.5 of \cite{buehlmannvdGeer2011}. In particular,
the diagram on page 177
exhibits compatibility as the weakest condition that still allows
oracle bounds for
prediction and reconstruction by the LASSO for the $\ell^2$- and $\ell
^1$-norms. The results on posterior
contraction and model selection presented below are in the same spirit.
In addition we consider contraction with respect to the $\ell^\infty$-norm,
and for (only) the latter we employ the more restrictive mutual
coherence number,
analogously to the study of \cite{lou} of the LASSO and the Dantzig
estimator under the
supremum norm. Thus mutual coherence is useful in two ways: it may
provide a simple (albeit crude)
way to bound the other indices, and it may allow to use stronger norms.
Direct verification of compatibility may be preferable, as this applies
to a much broader set of
regression matrices.

The following well-studied examples may help appreciate the discussion:

\begin{exam}[(Sequence model)]
In the sequence model of Example~\ref{ExampleSequenceModel}
the regression matrix $X$ is the identity, and hence the compatibility numbers
are 1, and the mutual coherence number is zero. This is
the optimal situation, under which all results below are valid.
(The compatibility numbers are maximally 1,
as follows by evaluating them with a unit vector.)

Regression with orthogonal design
can be transformed to this situation.
\end{exam}

\begin{exam}[(Response model)] In the response setting
of Example~\ref{ExampleRegressionModel} it is reasonable to assume that
the entries of $X$ are i.i.d. random variables. Under exponential
moment conditions, it is shown in \cite{Cai}
that in this situation and for not extremely large $p$ the mutual
coherence number
is with high probability bounded by a multiple of $(n/\log p)^{-1/2}$.
[Specifically, this is true for $\log p=o(n)$ or $\log p=o(n^{\alpha
/(4+\alpha)})$ if the entries
are bounded or possess an exponential moment of order $\alpha$,
resp.] In view of Lemma~\ref{LemmaCoherence} the compatibility
and sparse singular value indices of
models up to dimension a multiple of $\sqrt{n/\log p}$ are then
bounded away from zero.
This implies that the results on model selection and $\ell^1$- and
$\ell^2$-contraction rates
in the following certainly apply if the number of nonzero regression
coefficients is smaller than this
order. For a survey on more recent results on lower bounds of the
compatibility number and the smallest sparse eigenvalue, see Section~6.2 of \cite{vandegeermuro2014}.
\end{exam}

\begin{exam}
By scaling the columns of the design matrix it can be ensured that the
$(p\times p)$-matrix
$C:=X^tX/n$ has unit diagonal. Then $\|X\|=\sqrt{n}$, and the
off-diagonal elements $C_{i,j}$ are the correlations
between the columns.

It is shown in \cite{ZhaoYu} that if $C_{i,j}$ is equal to a constant
$r$ with
$0<r<(1+c s)^{-1}$, or $|C_{i,j}|\le c/(2s-1)$, for every $i\neq j$,
then models up to dimension $s$
satisfy the ``strong irrepresentability condition'' and hence are
consistently estimable. Since these
examples satisfy the mutual coherence condition, up to a constant,
these examples are also
covered in the present paper, for every norm and aspect considered.

As another example, Zhao and Yu \cite{ZhaoYu} consider correlations
satisfying $C_{i,j} = \rho^{|i-j|}$, for
$0<\rho<1$ and $p=n$. In this case \emph{all} eigenvalues of $C$ are
bounded away
from zero by a margin that depends on $\rho$ only, whence the numbers
$\tilde\phi(s)$ are
bounded away from zero, for every dimension $s$. This implies that the
results on dimensionality,
model selection and $\ell^1$- and $\ell^2$-rates obtained below are
valid. On the other hand,
the mutual coherence number is equal to $\rho$, which excludes the
$\ell^\infty$-results.

As a final example, the authors of \cite{ZhaoYu} consider matrices $C$
that vanish except in small blocks
along the diagonal. Such matrices can also not be handled in general
through the mutual coherence number,
but do cooperate with the other sparse invertibility indices.
\end{exam}

\subsection{Dimensionality, contraction, model selection}
For simplicity the main results are stated in limit form, for $p, n\ra
\infty$.
More precise assertions, including precise values of ``large'' constants,
can easily be deduced from the proofs.

The results are obtained under the assumption of Gaussian noise in model
\eqref{Model}. In fact, as indicated in Remark~1
in the supplementary material \cite{Supplement},
many of the assertions are robust under misspecification of the error
distribution. 

The first theorem shows that the posterior distribution does not
overshoot the true dimension of the
parameter by more than a factor. In the interesting case that $\la\ll
\ola$, this factor can be
simplified to $1+M/A_4$ for any constant $M>2$ if the true parameter is
compatible. The constant
$A_4$ comes from condition \eqref{assump.on_the_dim_prior}. As a
consequence, $1+M/A_4$ can be
made arbitrarily close to one by choosing a suitable prior on the
dimension. (Although the
convergence to zero in this and the following theorems is uniform, it
can be read off from the proofs
that the speed of convergence deteriorates for very small $\lambda$.
Also only the dominating terms in the
dependence of the dimension or contraction rate are shown. Thus the
theorems as stated
are not suitable for optimization over $\lambda$.
In particular, it should not be concluded that the smallest possible
$\lambda$ is optimal.)

\begin{thmm}[(Dimension)]
\label{thmm.mod_sel}
If $\la$ satisfies \eqref{eq.lambda_cond}
and $\pi_p$ satisfies \eqref{assump.on_the_dim_prior}
then, with $s_0=|S_{\be^0}|$ and for any $M>2$,
\[
\sup_{\be^0} \E_{\be^0}\Pi \biggl(\be:
|S_\be|> s_{0}+ \frac{M}{A_4} \biggl( 1+
\frac{16}{\phi(S_{0})^2}\frac{\la}{ \overline\la} \biggr) s_{0} \Big\Given Y \biggr)\ra0.
\]
\end{thmm}

The theorem is a special case of Theorem~\ref{thmm.mod_sel_gen} in
Section~\ref{SectionProofs}.
As all our results, the theorem concerns the full posterior
distribution, not only a measure of its center.
However, it may be compared to similar results for point estimators,
as in Chapter~7 of \cite{buehlmannvdGeer2011}.


The second theorem concerns the ability of the posterior distribution
to recover the
true model from the data. It gives rates of contraction of the
posterior distribution both regarding
\emph{prediction error} $\|X\be-X \be^0\|_2$ and regarding the
parameter $\be$ relative to the
$\ell^1$- and $\ell^2$- and $\ell^\infty$-distances.
Besides on the dimensionality, the rate depends on compatibility. Set
%
\begin{eqnarray}\label{DefPsis}
\overline\psi(S)&=&\overline\phi \biggl( \biggl( 2+\frac{3}{A_4}+
\frac{33}{\phi(S)^2}\frac{\la}{ \overline\la} \biggr) |S| \biggr),
\nonumber
\\[-8pt]
\\[-8pt]
\nonumber
\widetilde\psi(S)&=&\widetilde\phi \biggl( \biggl( 2+\frac{3}{A_4}+
\frac{33}{\phi(S)^2}\frac{\la}{ \overline\la} \biggr) |S| \biggr).
\end{eqnarray}
In the interesting case that $\la\ll\ola$, these numbers are
asymptotically bounded below
by $\overline\phi((2+\frac{4}{A_4})|S_\be|)$ and $\widetilde\phi
((2+\frac{4}{A_4})|S_\be|)$ if $\phi(S_\be)$ is
bounded away from zero. Thus the following theorem gives rates of recovery
that are uniform in true vectors $\be$ such that $\phi(S_\be)$ and
$\overline\phi((2+\frac{4}{A_4})|S_\be|)$
or $\widetilde\phi((2+\frac{4}{A_4})|S_\be|)$ are bounded away
from zero. [Again
the theorem, even though uniform in $\lambda$ satisfying \eqref
{eq.lambda_cond}, is meant
to be read for a given sequence of $\lambda$.]

\begin{thmm}[(Recovery)]
\label{TheoremRecovery}
If $\la$ satisfies \eqref{eq.lambda_cond},
and $\pi_p$ satisfies \eqref{assump.on_the_dim_prior},
then for sufficiently large $M$, with $S_0=S_{\be^0}$,
\begin{eqnarray*}
\sup_{\be^0} \E_{\be^0}\Pi \biggl(\be:\bigl \|X\bigl(\be-
\be^0\bigr)\bigr\|_2 > \frac
{M}{\overline\psi(S_{0})} \frac{\sqrt{|S_{0}| \log p}}{\phi(S_{0})} \Big| Y
\biggr) &\ra&0,
\\
\sup_{\be^0} \E_{\be^0}\Pi \biggl(\be: \bigl\|\be-
\be^0\bigr\|_1> \frac{M}{\overline
\psi(S_{0})^2} \frac{|S_{0}| \sqrt{\log p}}{\|X\|\phi(S_{0})^2}\Big| Y
\biggr)& \ra&0,
\\
\sup_{\be^0} \E_{\be^0}\Pi \biggl(\be:\bigl \|\be-
\be^0\bigr\|_2> \frac{M}{\widetilde
\psi(S_{0})^2} \frac{ \sqrt{|S_{0}| \log p}}{\|X\|\phi(S_{0})}\Big| Y
\biggr) &\ra&0.
\end{eqnarray*}
Furthermore, for every $c_0>0$, any $d_0<c_0^2(1+2/A_4)^{-1}/8$, and
$s_n$ with $\la s_n\sqrt{\log p}/\|X\|\ra0$, for sufficiently large $M$,
\[
\mathop{\sup_{\be^0:\phi(S_0)\ge c_0, \widetilde
\psi(S_0)\ge c_0}}_{
{|S_0|\le s_n, |S_0|\le d_0 \mc(X)^{-1} }} \E_{\be^0} \Pi
\biggl( \be: \bigl\|\be-\be^0\bigr\|_\infty> M\frac{\sqrt
{\log p}}{\|X\|} \Big| Y
\biggr) \ra0.
\]
\end{thmm}


The first three assertions of the
theorem are consequences of the following theorem of oracle type,
upon choosing $\be^*=\be^0$ in this theorem. The fourth assertion is
proved in Section~\ref{SectionProofs} under the conditions of Theorem~\ref{thmm.BvM_type} below. In the framework of Example~\ref
{ExampleRegressionModel}, for instance say for i.i.d. Gaussian design
and $\la=1$, the fourth assertion is true with large probability
uniformly over sparse vectors such that $|S_0|\le s_n=o(\sqrt{n/\log{p}})$.

An \emph{oracle inequality} for the prediction error of a point estimator
$\widehat\be$ is an assertion that with large probability, and
some penalty function $\pen(\be)$,
\[
\bigl\|X\bigl(\widehat\be-\be^0\bigr)\bigr\|_2^2
\lesssim\inf_{\be^*} \bigl\|X\bigl(\be^*-\be ^0\bigr)
\bigr\|_2^2+\pen\bigl(\be^*\bigr);
\]
see, for example, \cite{buehlmannvdGeer2011}, Theorem~6.2, or \cite
{BRT09} for
the LASSO or the Dantzig selector.
Few oracle-type results for \emph{posterior measures} have been
developed. (The results of
\cite{babenkobelitser}, for projection estimators in white noise, are
close relatives.)
The following theorem is an example of such a statement.
Given compatibility it shows that the bulk of the vectors $\be$ in the
support of the
posterior distribution satisfy an oracle inequality with penalty $\pen
(\be)=|S_\be|$.

\begin{thmm}[(Recovery, oracle)]
\label{thmm.pred_and_l1}
If $\la$ satisfies \eqref{eq.lambda_cond},
and $\pi_p$ satisfies \eqref{assump.on_the_dim_prior},
then, for $\overline\psi$ and $\widetilde\psi$ given in \eqref{DefPsis},
there exists a constant $M$ 
such that uniformly over $\be^0$ and $\be^*$ with $|S_*|\leq|S_0|$,
where $S_*=S_{\be^*}$,
\begin{eqnarray*}
&&\E_{\be^0}\Pi \biggl(\be: \bigl\|X\bigl(\be-\be^0
\bigr)\bigr\|_2 > \frac{M}{\overline\psi(S_{0})} \biggl[\bigl\|X\bigl(\be^*-
\be^0\bigr)\bigr\| _2+\frac{\sqrt{|S_{*}|\log p}}{\phi(S_{*})} \biggr] \Big| Y \biggr)
\ra 0,
\\
&&\E_{\be^0}\Pi \biggl(\be:\bigl \|\be-\be^0
\bigr\|_1 > \bigl\|\be^*-\be^0\bigr\|_1\\
&&\hspace*{32pt}{}+\frac{M}{\overline\psi(S_{0})^2}
\biggl[\frac{\|X(\be^*-\be^0)\|_2^2}{\|X\|\sqrt{\log p}} +\frac{|S_{^*}|
\sqrt{\log p}}{\|X\|\phi(S_{*})^2} \biggr] \Big| Y \biggr)\ra 0,
\\
&&\E_{\be^0}\Pi \biggl(\be: \bigl\|\be-\be^0
\bigr\|_2 > \frac{M}{\|X\|\widetilde\psi(S_0)^2} \biggl[\bigl\|X\bigl(\be^*-\be^0
\bigr)\bigr\| _2+\frac{\sqrt{|S_{*}|\log p}}{\phi(S_{*})} \biggr] \Big| Y
 \biggr)\ra 0.
\end{eqnarray*}
%
\end{thmm}

Besides the choice $\be^*=\be^0$, which yields the first three assertions
of Theorem~\ref{TheoremRecovery}, other choices of $\be^*$ also give
interesting
results. For instance, in the sequence model of Example~\ref
{ExampleSequenceModel},
the choice $\be^*=0$ gives that
\[
\sup_{\be^0} \E_{\be^0}\Pi \bigl(\be: \bigl\|\be-
\be^0\bigr\|_2 > M\bigl\|\be ^0\bigr\|_2 | Y
\bigr)\ra0.
\]
For $\|\be^0\|_2^2$ smaller than $|S_{\be^0}|\log p$, this improves on
Theorem~\ref{TheoremRecovery}, by quantifying the rate in the sizes
and not
only the number of nonzero coordinates in $\be^0$.

The posterior distribution induces a
distribution on the set of models $S\subset\{1,2,\ldots, p\}$, which
updates the prior masses given to these models by \eqref{DefPrior}. It is
desirable that this puts most of its mass on the true model $S_{\be
^0}$. As the support of a
vector $\be^0$ is defined only in a qualitative manner by its
coordinates $\be_i^0$
being zero or not, this will not
be true in general. However, the following theorem shows, under (only
strong) compatibility, that the posterior distribution will not charge models
that are strict supersets of the true model, no matter the magnitudes
of the
nonzero coordinates in $\be^0$. This may be considered the effect
of model selection through the prior $\pi_p$, which under our assumptions
prefers smaller models, enough so that it will not add unnecessary coordinates
when all truly nonzero coordinates are present.

\begin{thmm}
[(Selection: no supersets)]
\label{TheoremSelectionNoSupersets}
If $\la$ satisfies \eqref{eq.lambda_cond},
and $\pi_p$ satisfies~\eqref{assump.on_the_dim_prior} with $A_4>1$,
then for every $c_0>0$ and any $s_n\leq p^a$ with\break $s_n\la\sqrt{\log
p}/ \|X\|\ra0$ and $a<A_4-1$,
\[
\mathop{\sup_{\be^0:\phi(S_0)\ge c_0}}_{
|S_0|\le s_n, \widetilde\psi(S_0)\ge c_0} \E_{\be^0}\Pi (
\be:S_\be\supset S_{\be^0}, S_\be\neq
S_{\be^0} | Y ) \ra0.
\]
\end{thmm}

A nonzero coordinate of $\be^0$ that is too close to zero cannot be
detected as being
nonzero by any method. Consequently, the posterior distribution may
well charge
models $S$ that contain only a subset of the true model $S_{\be^0}$
and possibly other
coordinates, which is not
excluded by the preceding theorem. The following theorem gives thresholds
for detection, which become smaller as the compatibility conditions
become stronger.
The theorem may be compared to results in terms of \emph{beta-min conditions}
for point estimators; see, for example, \cite{buehlmannvdGeer2011},
Corollary~7.6.

\begin{thmm}
[(Selection)]
\label{TheoremSelection}
If $\la$ satisfies \eqref{eq.lambda_cond},
and $\pi_p$ satisfies \eqref{assump.on_the_dim_prior},
then, for sufficiently large $M$,
\begin{eqnarray*}
\inf_{\be^0} \E_{\be^0}\Pi \biggl(\be:
S_\be\supset \biggl\{i: \bigl|\be_i^0\bigr|\ge
\frac{M}{\overline\psi(S_{0})^2} \frac{|S_{0}|
 \sqrt{\log p}}{\|X\|\phi(S_{0})^2} \biggr\} \Big| Y \biggr) &\ra&1,
\\
\inf_{\be^0} \E_{\be^0}\Pi \biggl(\be:
S_\be\supset \biggl\{i: \bigl|\be_i^0\bigr|\ge
\frac{M}{\widetilde\psi(S_{0})^2}
\frac{ \sqrt{|S_{0}| \log p}}{\|X\|\phi(S_{0})} \biggr\}\Big | Y \biggr) &\ra&1.
\end{eqnarray*}
Furthermore, for every $c_0>0$, any $d_0\le c_0^2(1+2/A_4)^{-1}/8$, and
any $s_n$ with $\la s_n\sqrt{\log p}/\|X\|\ra0$,
\[
\mathop{\inf_{\be^0:\phi(S_0)\ge c_0, \widetilde
\psi(S_0)\ge c_0}}_{
{|S_0|\le s_n, |S_0|\leq d_0 \mc(X)^{-1} }} \E_{\be^0} \Pi
\biggl( \be: S_\be\supset \biggl\{i: \bigl|\be_i^0\bigr|
\ge \frac{M \sqrt{\log p}}{\|X\|} \biggr\} \Big| Y \biggr)\ra1.
\]
\end{thmm}

By combining Theorems~\ref{TheoremSelectionNoSupersets} and~\ref
{TheoremSelection}
we see that under the assumptions of the theorems
the posterior distribution \emph{consistently selects} the correct model
if \emph{all} nonzero coordinates
of $\be^0$ are bounded away from 0 by the thresholds given in
Theorem~\ref{TheoremSelection}.
For $M$ as in the preceding theorem, let
\[
\widetilde{\Beta}= \biggl\{\be: \min_{i\in S_\be}|
\be_i|\ge\frac
{M}{\widetilde\psi(S)^2} \frac{ \sqrt{|S_\be| \log p}}{\|X\|\phi(S_\be)} \biggr\}.
\]
Define $\overline{ \Beta}$ similarly with $\sqrt{|S_\be|\log p}$ in
the threshold
replaced by $|S_\be|\sqrt{\log p}$ and with $\overline\psi$ instead
of $\widetilde\psi$.

\begin{cor}[(Consistent model selection)]
\label{cor.consistent_mod_selection}
If $\la$ satisfies \eqref{eq.lambda_cond},
and $\pi_p$ satisfies \eqref{assump.on_the_dim_prior} with $A_4>1$,
and $s_n\leq p^a$ such that $a<A_4-1$ and\break $s_n\la\sqrt{\log p}/\|X\|
\ra0$, then, for every $c_0>0$,
\[
\mathop{\inf_{\be^0\in\widetilde\Beta:\phi
(S_0)\ge c_0}}_{|S_0|\le s_n,
\widetilde\psi(S_0)\ge c_0}\E_{\be^0}\Pi (
\be: S_\be =S_{\be^0}| Y )\ra1.
\]
The same is true with $\widetilde\Beta$ and $\widetilde\phi$
replaced by
$\overline{\Beta}$ and $\overline\phi$.
\end{cor}

Consistent posterior model selection implies in particular, that the
model with the largest posterior mass is model selection consistent in
the frequentist sense. This can be established as in the proof of
Theorem~2.5 in \cite{ggv}.

\subsection{Distributional approximation}
In this section we show that the posterior distribution can be
approximated by a mixture
of normal distributions. Moreover, given consistent selection of the
true model,
this mixture collapses to a single normal distribution. We restrict to
what we shall refer to as the \textit{small lambda regime},
%
\begin{equation}
 \frac{\la}{\|X\|} |S_{\be^0}|\sqrt{\log p}\ra0.
\label{eq.cond_on_lambda_BvM}
\end{equation}
%
In this case the centering of the normal distributions does not depend
on the size of
scaling parameters $\la$. In contrast, in the ``large lambda regime,''
which includes
the usual order of magnitude of the smoothing parameter in the LASSO,
the posterior distribution mimics
the LASSO, and gives a biased reconstruction of the true parameter; see
Theorem~1 in the supplementary material \cite{Supplement}.

The small lambda regime includes a variety of possible choices
within our general assumption \eqref{eq.lambda_cond}.
A smaller value of $\la$ corresponds to a noninformative prior
on the nonzero coordinates of the parameter vector.
Here ``small'' is relative, depending on the model and the number
of observations.

\begin{exam}
[(Small lambda regime)]
For the minimal choice $\la=\|X\|/p$ in~\eqref{eq.lambda_cond} the
small lambda regime
\eqref{eq.cond_on_lambda_BvM} simplifies to $|S_{\be^0}|\ll p/\sqrt
{\log p}$.
Thus the regime applies to a wide range of true parameters.

In the sequence model with multiple observations given in
Example~\ref{ExampleSequenceModelMultiple} and the response
model of Example~\ref{ExampleRegressionModel}, we have $\|X\|=\si
_n^{-1}$ and
$\|X\|\sim n^{1/2}$, respectively, and $\la$ is in the small lambda regime
if $\la|S_{\be^0}|$ is much smaller than $1/(\si_n\sqrt{\log p})$
and $\sqrt{n/\log p}$,
respectively. The second allows $\la=O(1)$ if $|S_{\be^0}|\sqrt{\log
p/n}\ra0$.
\end{exam}

For a given model $S\subset\{1,\ldots, p\}$ let $X_S$ be the $n\times
|S|$-submatrix of the regression
matrix $X$ consisting of the columns $X_{\cdot,i}$ with $i\in S$,
and let $\widehat\be_{(S)}$ be a least square estimator in the
restricted model
$Y=X_S\be_S+\ep$, that is,
\[
\widehat\be_{(S)}\in\argmin_{\be_S\in\RR^S} \|Y-X_S
\be_S\|_2^2.
\]
In case the restricted model would be correctly specified, the least
squares estimator
would possess a $\mathcal{N}(\be_S^0,(X_S^tX_S)^{-1})$-distribution,
and the posterior distribution
(in a setting where the data washes out the prior) would be
asymptotically equivalent to a
$\mathcal{N}(\widehat\be_{(S)},(X_S^tX_S)^{-1})$-distribution, by
the Bernstein--von Mises theorem.
In our present situation, the posterior distribution is approximated by a
random mixture of these normal distributions, of the form
\[
\Pi^\infty(\cdot| Y) = \sum_{S \in\mathcal{S}_0} \widehat
w_S \mathcal{N}\bigl(\widehat\be _{(S)},
\bigl(X_S^tX_S\bigr)^{-1}\bigr)
\otimes\delta_{S^c}, 
\]
where $\delta_{S^c}$ denotes the Dirac measure at $0\in\RR^{S^c}$,
the weights $(\widehat w_S)_S$ satisfy
%
\begin{equation}
\label{EqDefWeightsw} \widehat w_S\propto \frac{\pi_p(s)}{{p\choose s}}
 \biggl(
\frac{\la}2 \biggr)^s (2\pi )^{s/2}
\bigl|X_S^tX_S\bigr|^{-1/2}
e^{({1}/2) \|X_S\widehat\be_{(S)}\|_2^2}1_{S\in\mathcal{S}_0}
\end{equation}
and, for a sufficiently large $M$
\[
\mathcal{S}_{0}= \biggl\{S: |S|\le\biggl(2+\frac{4}{A_4}
\biggr)|S_{\be^0}|,\bigl \| \be_{S^c}^0\bigr\|_1
\le M |S_{\be^0}|\sqrt{\log p}/\|X\| \biggr\}.
\]
The weights $(\widehat w_S)$ are a data-dependent probability
distribution on the
collection of models $\mathcal{S}_{0}$. The latter collection
can be considered a ``neighborhood'' of the support
of the true parameter, both in terms of dimensionality and the (lack
of) extension of the true parameter
outside these models.

A different way of writing the approximation $\Pi^\infty$ is
%
\begin{equation}
\Pi^\infty(B| Y) = \frac{\sum_{S\in\mathcal{S}_0}
({\pi_p(s)}/{{p\choose s}})
 ({\la}/2 )^s
\int_{B_S} e^{-({1}/2) \|Y-X_S\be_S\|_2^2} \,d\be_S}{
\sum_{S \in\mathcal{S}_0} ({\pi_p(s)}/{{p\choose s}})
({\la}/2 )^s
\int e^{-({1}/2) \|Y-X_S\be_S\|_2^2} \,d\be_S}, \label{eq.post_infty_2nd_repr}
\end{equation}
where $B_S=\{\be_S: (\be_S,0_{S^c})\in B\}$ is the intersection
(and not projection) of $B\subset\RR^p$ with the subspace $\RR^S$.
To see this, decompose $Y-X_S\be_S=(Y-X_S\widehat\be
_{(S)})+X_S(\widehat\be_{(S)}-\be_S)$, and observe that the two
summands are orthogonal. The Lebesgue integral $d\be_S$ can be interpreted
as an improper prior on the parameter $\be_S$ of model $S$, and the
expression as
a mixture of the corresponding posterior distributions, with model weights
proportional to the prior weights times $(\la/2)^s(2\pi)^{s/2}\int
e^{-({1}/2) \|Y-X_S\be_S\|_2^2} \,d\be_S$.
It follows that the Laplace priors $g_S$ on the nonzero coordinates
wash out
from the components of the posterior. On the other hand, they are still
visible in the
weights through the factors $(\la/2)^s$. In general,
this influence is mild in the sense that these factors will not change
the relative
weights of the models much.

\begin{thmm}[(Bernstein--von Mises, small lambda regime)]
\label{thmm.BvM_type}
If $\la$ satisfies \eqref{eq.lambda_cond}, 
and $\pi_p$ satisfies \eqref{assump.on_the_dim_prior}, then
for every $c_0>0$ and any $s_n$ with $s_n\la\sqrt{\log p}/\|X\|\ra0$,
\[
\mathop{\sup_{\be^0:\phi(S_0)\ge c_0}}_{
|S_0|\le s_n,\overline\psi(S_0)\ge c_0} \E_{\be^0}\bigl \|\Pi(
\cdot| Y)-\Pi^\infty(\cdot| Y) \bigr\| _{\TV}\ra0.
\]
\end{thmm}


\begin{cor}[(Limit under strong model selection)]
\label{cor.strong_mod_selection}
Under the combined assumptions of Corollary~\ref{cor.consistent_mod_selection}
and Theorem~\ref{thmm.BvM_type},
\[
\mathop{\sup_{\be^0\in\widetilde\Beta:\phi
(S_0)\ge c_0}}_{{|S_0|\le s_n,
\widetilde\psi(S_0)\ge c_0}} \E_{\be^0} \bigl\|\Pi(
\cdot| Y)- \mathcal{N} \bigl(\widehat\be_{(S_0)}, \bigl(X_{S_0}^tX_{S_0}
\bigr)^{-1} \bigr) \otimes\delta_{S_0^c} \bigr\|_{\TV}\ra0.
\]
\end{cor}


The distributional results imply that the spread in the posterior distribution
gives a correct (conservative) quantification of remaining uncertainty
on the parameter.
One way of making this precise is in terms of \emph{credible sets} for
the individual parameters $\be_j$. The marginal posterior distribution
of $\be_j$
is a mixture
$\hat\pi_j\delta_0+\hat H_j$ of a point mass at zero and a
continuous component $\hat H_j$.
Thus a reasonable \emph{upper 0.975 credible limit} for $\be_j$ is
equal to
\[
\hat R_j=\cases{ \hat H_j^{-1}(0.975),&\quad$
\mbox{if } 0.975\le\hat H_j(0)$,\vspace *{2pt}
\cr
0,& \quad$\mbox{if }
\hat H_j(0)\le0.975\le\hat H_j(0)+\hat\pi
_j$,\vspace*{2pt}
\cr
\hat H_j^{-1}(0.975-\hat
\pi_j),&\quad $\mbox{if } \hat H_j(0)+\hat\pi
_j \le0.975$.}
\]
It is not difficult to see that under the conditions of Corollary~\ref
{cor.strong_mod_selection},
$\P_{\be^0} ( \be^0_j \leq\hat R_j)\ra0.975$ if $j\in S_0$
and $\P_{\be^0} ( \be^0_j=0)\ra1$ if $j\notin S_0$.


\section{The LASSO is not fully Bayesian}
\label{sec.spars_and_LASSO}

The LASSO (cf. \cite{tib})
\begin{equation}
\label{EqLASSO} \hat{\be}^{\mathrm{LASSO}}_\la= \argmin_{\be\in\RR^p}
\bigl[ \|Y-X\be\| _2^2 + 2\la\|\be\|_1
\bigr]
\end{equation}
is the posterior mode for the prior that models the coordinates $\be
_i$ as an i.i.d. sample from a Laplace
distribution with scale parameter $\la$, and thus also possesses a
Bayesian flavor.
It is well known to have many desirable properties: it is computationally
tractable; with appropriately tuned smoothing parameter $\la$ it attains
good reconstruction rates; it automatically leads to sparse
solutions; by small adaptations it can be made consistent for model
selection under standard conditions.
However, as a Bayesian object it has a deficit: in the sparse setup the
full posterior distribution corresponding
to the LASSO prior does not contract at the same speed as its mode.
Therefore the full posterior distribution is useless for uncertainty
quantification,
the central idea of Bayesian inference.

We prove this in the following theorem, which we restrict to the
sequence model of Example~\ref{ExampleSequenceModel}, that is, model
\eqref{Model} with $X=I$ the identity matrix.
In this setting the LASSO estimator is known to attain the (near)
minimax rate
$s \log n$ for the square Euclidean loss over the
``nearly black bodies'' $\{\be: |S_\be|\le s\}$,
and a near minimax rate over many other sparsity classes as well,
if the regularity parameter $\la$ is chosen of the order $\sqrt{2\log{n}}$.
The next theorem shows that for this choice the LASSO posterior distribution
$\Pi_\la^{\mathrm{LASSO}}(\cdot| Y)$ puts
no mass on balls of radius of the order $\sqrt n/(\log n)^{1/2}$, which is
substantially bigger than the minimax rate $(s \log n)^{1/2}$ (except
for extremely dense signals).

Intuitively, this is explained by the fact that the
parameter $\la$ in the Laplace prior must
be large in order to shrink coefficients $\be_i$ to zero, but at the
same time reasonable
so that the Laplace prior can model the nonzero coordinates.
That these conflicting demands do not affect the good behavior of the LASSO
estimators must be due to the special geometric, sparsity-inducing form
of the posterior mode,
not to the Bayesian connection.

\begin{thmm} \label{lemlb}
Assume that we are in the setting of Example~\ref
{ExampleSequenceModel}. For any $\la=\la_n$ such that $\rn/\la_n\to
\infty$, there exists $\delta>0$ such that, as $n\to\infty$,
\[
\E_{\be^0=0} \Pi_{\la_n}^{\mathrm{LASSO}} \biggl(\be: \|\be
\|_2\le\delta \sqrt{n} \biggl(\frac{1}{\la_n}\wedge1 \biggr) \Big| Y
\biggr) \rightarrow0.
\]
\end{thmm}


\section{Prediction for arbitrary design}
\label{SectionArbitraryDesign}
The vector $X\be$ is the mean vector of the observation $Y$ in \eqref
{Model}, and one might guess
that this is estimable without identifiability conditions on the
regression matrix $X$. In this
section we show that the posterior distribution based on the prior \eqref
{DefPrior} can indeed solve
this \emph{prediction problem} at (nearly) optimal rates under no
condition on the
design matrix $X$. These results are inspired by \cite{DT07} and
Theorem~\ref{thmm-orapred} below can be seen as a full Bayesian
version of the
results on the PAC-Bayesian point estimators in the latter paper; see also
\cite{rigollet-tsybakov11} for prediction results for mixtures of
least-squares estimators.

We are still interested in the sparse setting, and hence the regression matrix
$X$ still intervenes by modeling the unknown mean vector $\E Y$ as a linear
combination of a small set of its columns.

First, we consider the case of priors \eqref{DefPrior} that
model the mean vector indirectly by modeling the set of columns and
the coefficients of the linear combination. The prior $\pi_p(s)$ comes
in through the constant
%
\begin{equation}
\label{predcpi} C_\pi= \sum_{s=0}^p
9^s\pmatrix{p \cr s}^{1/2}\sqrt{\pi_p(s)}.
%
\end{equation}
%
For the choice of prior on coordinates $\be_i$, the best results are
obtained with heavy-tailed densities $g$.
In general the rate depends on the Kullback--Leibler divergence
between the measure with distribution function $G_{S_0}$ (corresponding
to the prior density $g_{S_0}$) and the same measure shifted by $\be^0_{S_0}$.
Let $\mathrm{KL}$ be the Kullback--Leibler divergence, and set
%
\begin{equation}
\label{def.db} D_{\be^0} =\frac{{p \choose s_0}}{\pi_p(s_0)} e^{\mathrm{KL}(G_{S_0}(\cdot-\be^0_{S_0}),G_{S_0})
+( {1}/{2})\int\|X\be_{S_0}\|_2^2 \,dG_{S_0}(\be_{S_0}) }.
\end{equation}

\begin{thmm} \label{thmm-orapred}
For any prior $\pi_p$ and $C_\pi$ as in \eqref{predcpi},
any density $g$ that is symmetric about $0$,
any $\be^0,\be^*\in\RR^p$ and $r\ge1$,
\[
\E_{\be^0} \Pi \Bigl( \bigl\|X\bigl(\be-\be^0\bigr)
\bigr\|_2 > 7\bigl\|X\bigl(\be^*-\be^0\bigr)\bigr\|_2 + 4
\sqrt{\log\bigl(C_\pi^2 D_{\be^*}\bigr)}
+ 8\sqrt{r} | Y \Bigr) \leqa e^{-r}.
\]
\end{thmm}

%
%



If the prior on the dimension satisfies \eqref
{assump.on_the_dim_prior} with $A_4>1$, then $C_\pi$ is bounded in~$p$, and the rate for squared error loss is determined by
\begin{eqnarray*}
\rho_n\bigl(\be^0\bigr) &:=& \log D_{\be^0}
\\
&\lesssim&|S_{\be^0}|\log{p}+ \mathrm{KL}\bigl(G_{S_0}\bigl(\cdot-
\be^0_{S_0}\bigr),G_{S_0}\bigr) +
\frac{1}{2}\int\|X\be_{S_0}\|_2^2
\,dG_{S_0}(\be_{S_0}).
\end{eqnarray*}
This rate might be dominated by the Kullback--Leibler divergence for
large signal~$\beta^0$. However, for heavy tailed priors
$g$ the induced constraints on the signal to achieve the good rate
$|S_{\be^0}|\log{p}$
are quite mild.
Consider the prior distribution \eqref{DefPrior} with $g_S$ a product of
$|S|$ univariate densities $g$ of the form
%
\begin{equation}
\label{def.ht} g(x) \propto\frac{\la}{1+|\la x|^{\mu}},\qquad  x\in\RR, \la>0, \mu>3.
\end{equation}

\begin{cor} \label{corpredht}
If $\pi_p$ satisfies \eqref{eq.classofdimpriors} with $a\ge1$, and
$g$ is of the form \eqref{def.ht}
with $\la=\|X\|$ and $\mu>3$, then for sufficiently large $M$,
\[
\sup_{\be^0} \E_{\be^0}\Pi \bigl(\be\in
\RR^p: \bigl\|X\be-X\be^0\bigr\| _2^2> M
\rho_n\bigl(\be^0\bigr) | Y \bigr)\ra0,
\]
for $\rho_n(\be) = |S_\be|\log{p}  \vee\sum_{i\in S_\be}\log
 (1+\|X\|^\mu|\be_{i}|^\mu )$.
\end{cor}


\begin{rem}
The constant $7$ in Theorem~\ref{thmm-orapred} can be improved to
$4+\delta$, for an arbitrary $\delta>0$, by a slight adaptation of
the argument. Using
PAC-Bayesian techniques Dalalyan and Tsybakov \cite{DT07} obtain an
oracle inequality with leading constant $1$ for a
so-called pseudo-posterior mean: the likelihood in \eqref{bayes} is
raised to some power, which
amounts to replacing the $1/2$ factor by $1/\beta$. The ``inverse
temperature'' $\beta$ must
be taken large enough; the case $\beta=2$ corresponding to the Bayes
posterior as considered
here is not included; see also \cite{leungbarron06}.
\end{rem}

Theorem~\ref{thmm-orapred} and its corollary address the question of
achieving prediction with
no condition on $X$, and the same rate is achieved as in Section~\ref{sec.recovery} with the same type of priors, up to some slight loss
incurred only for true vectors $\be^0$ with very large entries. As
shown in the corollary, this slight dependence on $\be^0$ can be made
milder with flatter priors.
We now consider a different approach specifically targeted at the
prediction problem and which enables to remove dependency on the size
of the coordinates of $\be^0$ completely.

Because the prediction problem is concerned
only with the mean vector, and the columns of $X$ will typically be
linearly dependent,
it is natural to define the prior distribution directly on the corresponding
subspaces.
For any $S\subset\{1,\ldots,p\}$, let $\cY_S:=\{X\be, S_{\be
}\subseteq S\}$ be the subspace of $\RR^n$ generated by the columns
$X^j, j\in S$ of~$X$. Let $\cV$ denote the collection of all
\emph{distinct} subspaces $\cY_S$.

Define a (improper) prior $\Xi$ on $\RR^n$ by first selecting an
integer $t$ in $\{0,1,\ldots,n\}$ according to a prior $\pi_n$, next
given $t$ selecting a subspace $V\in\cV$ of dimension $t$ uniformly
at random among subspaces in $\cV$ of dimension $t$; finally, let $\Xi
$ given $V$ be defined as Lebesgue measure on $V$ if $\mbox{dim}(V)\ge
1$, and let $\Xi$ be the Dirac mass at $\{0\}$ for $V=\{0\}$. Note
that the posterior distribution $\Xi[\cdot| Y]$ is a well-defined
probability measure on $\RR^n$.

We choose, for a fixed $d \ge4$ (the numerical constant $4$ is for simplicity),
\begin{equation}
\label{def-pn} \pi_n(t) := \pi_{n,p}(t)=
\frac{e^{-dt\log{p}}}{ \sum_{t=0}^n
e^{-dt\log{p}}},\qquad t=0,1,\ldots,n.
\end{equation}
Let $V^0:=\cY_{S_{\be^0}}$ and $t_0$ be the dimension of $V^0$.

\begin{thmm} \label{thmm-improper}
Let $\Xi$ be the improper prior on $\RR^n$ defined above with $\pi
_n$ as in~\eqref{def-pn}. For $M$ large enough,
\[
\sup_{\be^0} \E_{\be^0}\Xi\bigl[ \gamma\in
\RR^n, \bigl\|\gamma-X\be ^0\bigr\|_2^2 >
M(t_0\vee1)\log{p} | Y\bigr] \to0.
\]
\end{thmm}
The result is uniform in $\be^0\in\RR^p$. Also, note that $t_0\le
|S_{\be^0}|=s_0$ and that one may have $t_0 =o(s_0)$. The obtained
rate thus may improve on the previous prediction rates.
It has a simple interpretation: up to an additional logarithmic factor,
it is the rate of the natural estimate $\gamma^*=\mbox{Proj}_{V^0}Y$
if the true subspace $V^0$ is known, where $\mbox{Proj}_{V^0}$ denotes
the orthogonal projection in $\RR^n$ into the subspace $V^0$.

\section{Computational algorithms}
\label{SectionComputationalAspects}
In this section we survey computational methods to compute posterior
distributions
in the regression model \eqref{Model} based on model selection priors
\eqref{DefPrior}. In most cases, this is a ``spike and slab'' prior,
as discussed
in Example~\ref{example.slabspike}, implemented with auxiliary 0--1
variables that indicate
whether a parameter $\be_j$ is included in the model or not.
The slab distribution is typically chosen a scale mixture of Gaussian
distributions, which
may include the Laplace law, which is an exponential mixture. Most
implementations
also allow an unknown error variance (which is taken to be unity in the
present paper),
with the inverse gamma distribution as the favorite prior.

For low-dimensional regression problems, computation of the posterior given
mixture priors was studied by many authors, including \cite
{mitchellbeauchamp,georgefoster,ishwaranrao,Dellaportas,yuanglin,ScottBerger10,JiSchnmidler}.
Higher-dimensional settings have been
considered recently: most of the following papers have appeared in the
last five years,
and a number of them are preprints.

Several authors
\cite
{HansDobraWest,StingoVannucci,LiZhang,BottoloRichardson,BottoloRichardson2,RichardsonBottolo}
have implemented \emph{MCMC} schemes to simulate from the posterior
distribution,
coupled with \emph{stochastic search} algorithms that limit the model
space, so as
to alleviate the curse of dimensionality.
Besides computation time, monitoring the convergence of the samplers is
an issue.
For higher dimensions it is impossible to sample from the complete
model space,
but this should also not be necessary, as in sparse situations
the posterior will concentrate on lower-dimensional
spaces, as is also apparent from our theoretical results.
Bottolo et al. \cite{BottoloRichardson2} provide ready-made
software, which runs on
dimensions up to several thousands. The same authors have
also exploited hardware solutions, such as graphical processing units,
to speed up computations in genomic data analyses.

\emph{Sequential Monte Carlo methods} or \emph{particle filters} can
be viewed as MCMC schemes
that can more readily
incorporate correct moves in the model space that ensure good approximation
to the posterior distribution. In \cite{ShiDunson,SchaferChopin} such methods
are shown to perform well for model selection in regression models with
up to
hundreds of covariates.

The \emph{shrinkage-thresholding Metropolis adjusted Langevin
algorithm} (or STMALA)
introduced in \cite{2013arXiv1312.5658S} is another variation on
earlier MCMC algorithms,
targeted to work for $p>n$, in, for instance, imaging applications.
It jointly samples a model and a regression vector in this model,
using proposals based on the gradient of the logarithm of the smooth
part of the
posterior distribution (as in MALA) combined with applying a
shrinkage-thresholding operator to set
coordinates to zero. Geometric convergence of the algorithm,
which is capable of moving between rather distant models, is guaranteed
for slab prior densities of the form $\propto\exp(-\lambda\|\beta\|
_1-\mu\|\beta\|_2^2)$, where $\mu>0$.
Illustrations showing good practical performance are given in \cite
{2013arXiv1312.5658S} (Section~5.2)
for values of $(n,p)$ equal to $(100,200)$ or $(39, 300)$.

An alternative to simulation from the exact posterior is to compute an
exact, analytic
approximation to the posterior. A relatively simple and computationally
efficient
\emph{variational Bayes approximation} is proposed in \cite{Ormerod}
and is shown
to perform satisfactorily, but examples in the paper are limited to
cases where $p\le n$.

By relaxing the spike at zero to a Gaussian distribution with small
variance, Ro{\v{c}}kov{\'a} and George
\cite{RockovaGeorge} succeeded in reducing computations of aspects of
the posterior
distribution, such as means and moments, to iterations of an efficient
\emph{EM-algorithm}.
They show good performance with exponentially decreasing priors on
model dimension, as considered
in the present paper.

Closely related to the spike and slab prior is \emph{exponential
weighting}, where each of the $2^p$ models
is given a prior weight, which is then updated with the likelihood
function. A survey and
numerical simulations in high-dimensional settings using the
Metropolis--Hastings algorithm can be
found in \cite{rigollet2012}. Stable reconstructions in dimensions
up to $p=500$, $n=200$ and sparsity level $s_0=20$ are shown to require
usually no more than
2000 iterations.

An (empirical, pseudo-) Bayes approach with a spike and Gaussian slabs centered
at the least square solutions of the underlying model is implemented
in \cite{2014arXiv1406.7718M}. The algorithm, which can be initialized
at the LASSO estimator,
is shown to perform well for $n$ up to 100 and $p$ up to 1000.
Because the slabs are centered on data-based quantities, the target of
this algorithm
is different from the posterior distribution in the present paper.
However, since the prior puts mass on all models, its computational complexity
is comparable to the procedure in the present paper.

For the sequence model of Example~\ref{ExampleSequenceModel}, an algorithm
to compute posterior quantities such as modes and quantiles based on
\emph{generating polynomials}
is implemented in \cite{spa}. This is efficient in terms of
computation time, but requires large memory.
Up to $n=p=500$ standard software and hardware suffice. The method may be
extended to other designs by making suitable transformations \cite{Belitser}.

\section{Proofs for Section~\texorpdfstring{\protect\ref{sec.recovery}}{2}}
\label{SectionProofs}

Denote by $p_{n,\be}$ the density of the $\mathcal{N}(X\be,I)$-distribution,
and the corresponding log likelihood ratios by
%
\begin{equation}
\Lambda_{n,\be,\be^*}(Y) =\frac{p_{n,\be}}{p_{n,\be^*}}(Y)
=e^{-({1}/2) \|X(\be-\be^*)\|_2^2
+(Y-X\be^*)^t X(\be-\be^*)}.
\label{eq.lr_representation}
\end{equation}

\begin{lem}
\label{lem.expl_bd_for_denom}
For $p$ sufficiently large and any $\be^* \in\RR^p$, with support
$S_*$ and $s_*:=|S_*|$, and
$\Pi$ given by \eqref{DefPrior} with $g_S$ a product of Laplace
densities with scale~$\la$, we have, almost surely,
\[
\int\Lambda_{n,\be,\be^*}(Y) \,d\Pi(\be) \geq\frac{\pi_p(s_*)}{p^{2s_*}}e^{-\la\|\be^*\|_1}e^{-1}.
\]
\end{lem}

\begin{pf}
For $s_*=0$ the right-hand side is $\pi_p(0)e^{-1}$, while the
left-hand side is bounded
below by $\Lambda_{n,0,0}\pi_p(0)=\pi_p(0)$, by \eqref{DefPrior}. Thus
we may assume that $s_*\ge1$.

First we prove that for any set $S$ and $s=|S|>0$,
%
\begin{equation}
\int_{\|\be_S\|_1\leq r} g_S(\be_S) \,d
\be_S = e^{-\la r}\sum_{k=s}^\infty
\frac{(\la r)^k}{k!} \geq e^{-\la r} \frac{(\la r)^s}{s!}. \label{eq.int_explicit}
\end{equation}
If $(L_i)_{i=1,\ldots,s}$ are i.i.d. random variables
with the Laplace distribution with scale parameter $\la$,
then $(|L_i|)_{i=1,\ldots,s}$ are i.i.d. exponential variables of the
same scale.
Hence the left-hand side of the display, which is equal to\break $\P
(\sum_{i=1}^s |L_i|\leq r )$,
is the probability that the first $s$ events of a Poisson process of intensity
$\la$ occur before time $r$. This is identical to the probability that
the Poisson
process has $s$ or more events in $[0,r]$, which is the sum in the display.

By \eqref{DefPrior}, the left-hand side of the lemma
is bounded below by
\begin{eqnarray*}
&&\frac{\pi_p(s_*)}{{p\choose s_*}} \int\Lambda_{n,\be,\be^*}(Y) g_{S_*}(
\be_{S_*}) \,d\be_{S_*}
\\
&&\qquad \geq\frac{\pi_p(s_*)}{{p\choose s_*}}
e^{-\la\|\be^*\|_1} \int e^{-({1}/2) \|X b_{S_*}\|_2^2+(Y-X\be^*)^t X
b_{S_*}}g_{S_*}(b_{S_*})
\,d b_{S_*},
\end{eqnarray*}
by \eqref{eq.lr_representation},
the change of variables $\be_{S_*}-\be^*_{S_*}\ra b_{S_*}$ and the inequality
$g_{S_*}(\be_{S_*})\ge e^{-\la\|\be^*\|_1}g_{S_*}(b_{S_*})$.
The finite measure $\mu$ defined by the identity $d\mu=\break  \exp(-\frac{1}2 \|X b_{S_*}\|_2^2) g_{S_*}(b_{S_*}) \,db_{S_*}$ is symmetric
about zero, and hence the mean of $b_{S_*}$ 
relative to $\mu$ is zero. Let $\bar\mu$ denote the normalized
probability measure corresponding to $\mu$, that is, $\bar\mu:=\mu
/\mu(\RR^{|S_*|})$. Let $\mathbb{E}_{\bar\mu}$ denote the
expectation operator with respect to $\bar\mu$. Define
$Z(b_{S_*}):=(Y-X\be^*)^t X b_{S_*}$.
By Jensen's inequality
$\E_{\bar\mu} \exp(Z)\ge\exp(\E_{\bar\mu} Z)$. However, $\E
_{\bar\mu} Z=0$, by the just mentioned symmetry of~$\mu$. So the
last display is bounded below by
\[
\frac{\pi_p(s_*)}{{p\choose s_*}}e^{-\la\|\be^*\|_1} \int
 e^{-({1}/2) \|X b_{S_*}\|_2^2} g_{S_*}(b_{S_*})
\,db_{S_*},
\]
almost surely. Using that $\|X\be\|_2=\|\sum_{i=1}^p\be_iX_{\cdot,i}\|
_2\le\|\be\|_1\|X\|$, and
then \eqref{eq.int_explicit}, we find that the integral in the last
display is bounded below by
\[
e^{-1/2}\int_{\|X\| \| b_{S_*}\|_1\leq1} g_{S_*}(b_{S_*})
\,db_{S_*} \geq e^{-1/2}\frac{e^{-\la/\|X\|}  (\la/\|X\| )^{s_*}}{s_*!}.
\]
With \eqref{eq.lambda_cond}, $e^{-\la/\|X\|} (\la/\|X\|
)^{s_*}$ is bounded from below by $e^{-1/2} p^{-s_*}$, if $\lambda/\|
X\| \linebreak\leq1/2$ and by $e^{-4\sqrt{\log p}}2^{-s_*}$, if
$\lambda/\|X\|\geq1/2$. Since $s_*>0$ and $e^{-4\sqrt{\log p}}$
decays to zero slower than any polynomial power of $p$, we find
$e^{-\la/\|X\|} (\la/\|X\| )^{s_*} \linebreak\geq e^{-1/2}
p^{-s_*}$ in both cases, provided that $p$ is sufficiently large. The
lemma follows upon substituting these bounds and the bound
${p\choose s}s!\le p^s$ in the display.
\end{pf}

\begin{lem}
\label{lem.com}
For any $\be, \be^*\in\RR^p$ and random variable $0\le U=U(Y)\le1$,
\[
\E_{\be}U \leq e^{({1}/2) \|X(\be-\be^*)\|_2^2} (\E_{\be
^*}U
)^{1/2}.
\]
\end{lem}

\begin{pf}
Write the left-hand side as $\E_{\be^*} [U\Lambda_{n,\be,\be
^*} ]$, and use the Cauchy--Schwarz inequality; see, for example,
\cite{spa}, Lemma~6.1.
\end{pf}

\begin{lem}
\label{eq.mcal_T_def}
\[
\P_{\be^0} \bigl( \bigl\|X^t\bigl(Y-X\be^0\bigr)
\bigr\|_\infty> 2\sqrt {\log p} \|X\| \bigr)\le\frac{2}p. 
\]
\end{lem}

\begin{pf}
Under the probability measure $\P_{\be^0}$ the vector $\ep=Y-X\be^0$
possesses an $n$-dimensional standard normal distribution, whence
the $p$ coordinates of the vector $X^t\ep$ are normal with
variances $(X^tX)_{i,i}\le\|X\|$.
Now $\P (\|X^t\ep\|_\infty>x )\le\sum_{i=1}^p \P
((X^t\ep)_i>x )$,
which can be bounded by the tail bound for the normal distribution.
\end{pf}

\begin{thmm}[(Dimension, general result)]
\label{thmm.mod_sel_gen}
If $\la$ satisfies \eqref{eq.lambda_cond} and
the prior $\pi_p$ satisfies \eqref{assump.on_the_dim_prior}, then for
any $M>2$,
\begin{eqnarray*}
&&\sup_{\be^0, \be^*} \E_{\be^0}\Pi \biggl(\be:
|S_\be|> |S_{*}|+ \frac{M}{A_4} \biggl( 1+
\frac{16}{\phi(S_{*})^2}\frac{\la}{ \overline\la} \biggr) |S_{*}|+ \frac{M}{A_4}
\frac{\|X(\be^0-\be^*)\|_2^2}{\log p} \Big| Y \biggr)\\
&&\qquad\ra0.
\end{eqnarray*}
\end{thmm}

\begin{pf}
By the definition of $\overline\la$ in \eqref{eq.lambda_cond}
and Lemma~\ref{eq.mcal_T_def}, the complement of the event
$\cT_{0}=\{\|X^t(Y-X\be^0)\|_\infty\le \ola\}$ has $\P_{\be
^0}$-probability bounded by $2/p$.
By combining this with Lemma~\ref{lem.com} we see that for any $\be^*$
and any measurable set $B\subset\RR^p$,
%
\begin{equation}
\E_{\be^0}\Pi(B| Y) \leq e^{({1}/2) \|X(\be^0-\be^*)\|_2^2} \bigl(\E_{\be^*}
\bigl[\Pi(B| Y) \mathbb{I}_{\cT_{0}} \bigr] \bigr)^{1/2} +
\frac{2}p. \label{eq.main_decomp_mod_sel_gen}
\end{equation}
By Bayes's formula followed by Lemma~\ref{lem.expl_bd_for_denom},
with $\Lambda_{n,\be,\be^*}(Y)$ the likelihood ratio given in \eqref
{eq.lr_representation},
%
\begin{eqnarray}
\label{eq.post_exp_bd_gen} \Pi(B| Y) &=&\frac{\int_B \Lambda_{n,\be,\be^*}(Y) \,d\Pi(\be)}{
\int\Lambda_{n,\be,\be^*}(Y) \,d\Pi(\be)}
\nonumber
\\[-8pt]
\\[-8pt]
\nonumber
&\leq& \frac{e p^{2s_*}}{\pi_p(s_*)} e^{\la\|\be^*\|_1} \int_B
e^{-({1}/2) \|X(\be-\be^*)\|_2^2 +(Y-X\be^*)^t X(\be-\be
^*)} \,d\Pi(\be).
\end{eqnarray}
Using H\"older's inequality $|\al^t\be|\le\|\al\|_\infty\|\be\|
_1$ and
the Cauchy--Schwarz inequality, we see that on the event $\cT_{0}$,
%
\begin{eqnarray}
\label{EqEstimateInnerProduct} &&\bigl(Y-X\be^*\bigr)^t X\bigl(\be-\be^*\bigr)\nonumber\\
&&\qquad =
\bigl(Y-X\be^0\bigr)^tX\bigl(\be-\be^*\bigr)+\bigl(X
\be^0-X\be^*\bigr)^tX\bigl(\be-\be^*\bigr)
\\
&&\qquad\le\ola\bigl\|\be-\be^*\bigr\|_1+\bigl \|X\bigl(\be^0-\be^*\bigr)
\bigr\|_2^2 +\tfrac{1}{4} \bigl\|X\bigl(\be-\be^*\bigr)
\bigr\|_2^2=:L(\be).\nonumber
\end{eqnarray}
Therefore, on the event $\cT_0$, the expected value under
$\mathbb{E}_{\be^*}$ of the integrand on the right-hand side of
\eqref{eq.post_exp_bd_gen} is bounded above by
\begin{eqnarray*}
&&e^{-({1}/2) \|X(\be-\be^*)\|_2^2}
\E_{\be^*} \bigl[e^{(1-{\la}/{(2\ola)}) (Y-X\be^*)^tX(\be-\be
^*)}
\mathbb{I}_{\cT_{0}} \bigr] e^{({\la}/{(2\ola)}) L(\be)}
\\
&&\qquad =e^{-({1}/2 )(1-(1-{\la}/{(2\ola)})^2 ) \|X(\be-\be
^*)\|_2^2} e^{({\la}/{(2\ola)}) L(\be)}
\\
&&\qquad \le e^{\|X(\be^0-\be^*)\|_2^2}
e^{-({\la}/{(8\overline\la)}) \|X(\be-\be^*)\|_2^2+({\la}/2)
\|\be-\be^*\|_1},
\end{eqnarray*}
where we use that $\la\le2\ola$. It follows that the expected value
$\E_{\be^*} [\Pi(B| Y) \mathbb{I}_{\cT_0} ]$ under $\be^*$ of
\eqref{eq.post_exp_bd_gen} over $\cT_0$ is bounded above by
%
\begin{equation}\quad
\le \frac{e p^{2s_*}}{\pi_p(s_*)} e^{\la\|\be^*\|_1} \int_B
e^{\|X(\be^0-\be^*)\|_2^2}
e^{-({\la}/{(8\overline\la)}) \|X(\be-\be^*)\|_2^2+({\la}/2)
\|\be-\be^*\|_1} \,d\Pi(\be). \label{EqExpectedOverTo}\hspace*{-12pt}
\end{equation}
By the triangle inequality,
\begin{eqnarray*}
\bigl\|\be^*\bigr\|_1+\tfrac{1}2\bigl\|\be-\be^*\bigr\|_1 &\le&\|
\be_{S_*}\|_1+ \tfrac{3}2 \bigl\|\be_{S_*}-
\be^*\bigr\|_1+ \tfrac{1}2\| \be_{S_*^c}\|_1
\\
&\le&-\tfrac{1}4 \bigl\|\be-\be^*\bigr\|_1+ \|\be\|_1,
\end{eqnarray*}
for $7\|\be_{S_*}-\be^*\|_1\le\|\be_{S_*^c}\|_1$, as is seen by splitting
the norms on the right-hand side over $S_*$ and $S_*^c$.
If $\|\be_{S_*^c}\|_1< 7\|\be_{S_*}-\be^*\|_1$, then we write
$3/2=2-1/2$ and
use the definition of the compatibility number $\phi(S_*)$ to find that
\begin{eqnarray*}
&&\|\be_{S_*}\|_1+ \frac{3}2 \bigl\|
\be_{S_*}-\be^*\bigr\|_1+\frac{1}2\|\be
_{S_*^c}\|_1
\\
&&\qquad \le\|\be_{S_*}\|_1+ 2 \frac{\|X(\be-\be^*)\|_2s_*^{1/2}}{\|X\|
\phi(S^*)} -
\frac{1}2 \bigl\|\be_{S_*}-\be^*\bigr\|_1+\frac{1}2
\|\be_{S_*^c}\|_1
\\
&&\qquad \le\|\be_{S_*}\|_1+\frac{1}{8\ola}\bigl\|X\bigl(\be-\be^*
\bigr)\bigr\|_2^2+\frac
{8s_*\ola}{\|X\|^2\phi(S_*)^2} -\frac{1}4 \bigl\|\be-
\be^*\bigr\|_1+\|\be\|_1.
\end{eqnarray*}
We combine the last three displays to see that
\eqref{EqExpectedOverTo} is bounded above by
\[
\frac{e p^{2s_*}}{\pi_p(s_*)}e^{\|X(\be^0-\be^*)\|_2^2} e^{
{8\la\ola s_*}/{(\|X\|^2 \phi(S_*)^2)}} \int_B
e^{-({\la}/4)\|\be-\be^*\|_1+\la\|\be\|_1} \,d\Pi(\be). 
\]
For the set $B=\{\be: |S_\be|> R\}$ and $R\geq s_*$, the integral in
this expression is bounded above by
\begin{eqnarray*}
&&\sum_{S:|S|> R} \frac{\pi_p(s)}{{p\choose s}} \biggl(
\frac{\la}2 \biggr)^s \int e^{-({\la}/4) \|\be_S-\be^*\|
_1} \,d\be_S\\
&&\qquad \leq\sum_{s=R+1}^p \pi_p(s)
4^s
\\
&&\qquad\le\pi_p(s_*)4^{s_*} \biggl(\frac{4A_2}{p^{A_4}}
\biggr)^{R+1
-s_*}\sum_{j=0}^\infty
\biggl(\frac{4A_2}{p^{A_4}} \biggr)^j,
\end{eqnarray*}
by assumption \eqref{assump.on_the_dim_prior}.
Combining the preceding with \eqref{eq.main_decomp_mod_sel_gen}, we
see that
\[
\E_{\be^0}\Pi(B| Y)
\lesssim(4p)^{s_*} e^{\|X(\be^0-\be^*)\|_2^2
+{4\la\ola s_*}/{(\|X\|^2 \phi(S_*)^2)}}
\biggl(\frac{4A_2}{p^{A_4}} \biggr)^{(R+1 -s_*)/2}+ \frac{2}p.
\]
Using that $\ola^2=4\|X\|^2\log p$, we can infer the theorem by choosing
$R=s_*+M A_4^{-1}(\|X(\be^0-\be^*)\|_2^2/\log p+ s_*+16s_*(\la/\ola
)/\phi(S_*)^2)$ for fixed $M>2$.
\end{pf}

\begin{pf*}{Proof of Theorem~\ref{thmm.pred_and_l1}}
By Theorem~\ref{thmm.mod_sel_gen} the posterior distribution is asymptotically
supported on the event $E:=\{\be: |S_\be|\le D_* \wedge D_0\}$, for
%
\begin{equation}
\label{DefDimensionStar} D_*= \biggl( 1+ \frac{3}{A_4}+ \frac{33}{A_4\phi(S_{*})^2}
\frac
{\la}{ \overline\la} \biggr) s_{*}+ \frac{3}{A_4}\frac{\|X(\be
^0-\be^*)\|_2^2}{\log p}
\end{equation}
and $D_0$ the same expression with $\be^*$ replaced by $\be^0$. Thus
it suffices
to prove that the intersections of the events in the theorem with the
event $E$ tends to zero.
By combining \eqref{eq.post_exp_bd_gen}, \eqref
{EqEstimateInnerProduct} and
the inequality $\la\|\be^*\|_1\le2\ola\|\be-\be^*\|_1+\la\|\be\|
_1$, we see
that on the event $\cT_{0}=\{\|X^t(Y-X\be^0)\|_\infty\le \ola\}$,
the variable $\Pi(B| Y)$ is
bounded above by
\[
\frac{e p^{2s_*}}{\pi_p(s_*)} \int_B e^{-({1}/4) \|X(\be-\be^*)\|_2^2 +3\ola\|\be-\be^*\|_1+\|
X(\be^0-\be^*)\|_2^2+\la\|\be\|_1} \,d\Pi(\be).
\]
By Definition~\ref{assump.uniform_comp} of the uniform compatibility number,
%
\begin{eqnarray}\qquad
\label{EqInvokeCompCond} (4-1)\ola\bigl\|\be-\be^*\bigr\|_1 &\le&
\frac{4\ola\|X(\be-\be^*)\|_2|S_{\be-\be^*}|^{1/2}}{
\|X\|\overline\phi(|S_{\be-\be^*}|)} -\ola
\bigl\|\be-\be^*\bigr\|_1
\\
&\le&\frac{1}8 \bigl\|X\bigl(\be-\be^*\bigr)\bigr\|_2^2+
\frac{32\ola^2|S_{\be-\be^*}|}{
\|X\|^2\overline\phi(|S_{\be-\be^*}|)^2} -\ola\bigl\|\be-\be^*\bigr\|_1.
\end{eqnarray}
Since $|S_{\be-\be^*}|\le|S_\be|+s_*\le D_*\wedge D_0+s_*$, on the
event $E$ and $s_*\leq s_0$ by assumption, it follows from \eqref
{DefPsis} that
for a set $B\subset E$,
%
\begin{eqnarray}
\label{eq.post_exp_bd_gen2} \Pi(B| Y) \mathbb{I}_{\cT_0}&\le&\frac{e p^{2s_*}}
{\pi_p(s_*)}
e^{\|X(\be^0-\be^*)\|_2^2+{32\ola^2 (D_*+s_*)}/{(\|X\|
^2\overline\psi(S_0)^2)}}
\nonumber
\\[-8pt]
\\[-8pt]
\nonumber
& &{}\times\int_B e^{-({1}/8) \|X(\be-\be^*)\|_2^2 -\ola\|\be-\be
^*\|_1+\la\|\be\|_1} \,d\Pi(\be).
\end{eqnarray}
Since $\P_{\be^0}(\cT_0)\leq2/p$ it suffices to show that the
right-hand side tends to zero
for the relevant event $B$.

\emph{Proof of first assertion}.
On the set $B:=\{\be\in E: \|X(\be-\be^0)\|_2>4 \|X(\be^*-\be^0)\|
_2+R\}$, we have
$\|X(\be-\be^*)\|_2^2> 9\|X(\be^*-\be^0)\|_2^2+R^2$, by the
triangle inequality. Note that $\pi_p(s_*)\geq(A_1p^{-A_3})^{s_*}\pi
_p(0)$. It follows that for the set $B$, the preceding display is
bounded above by
\begin{eqnarray*}
&&\frac{e p^{2s_*}}{\pi_p(s_*)} e^{{32\ola^2 (D_*+s_*)}/
{(\|X\|^2\overline\psi(S_0)^2)}}e^{-({1}/8) R^2} \int e^{-\ola\|\be-\be^*\|_1+\la\|\be\|_1} \,d\Pi(
\be)
\\
&&\qquad \lesssim p^{(2+A_3)s_*} A_1^{-s_*}
 e^{{32\ola^2 (D_*+s_*)}/{(\|X\|^2\overline\psi(S_0)^2)}}e^{-({1}/8) R^2}
\sum_{s=0}^p\pi_p(s)2^s,
\end{eqnarray*}
by \eqref{assump.on_the_dim_prior} and
a calculation similar to the proof of Theorem~\ref{thmm.mod_sel_gen}.
For
\begin{eqnarray*}
\frac{1}8R^2&=&(3+A_3)s_*\log p +
\frac{32\ola^2 (D_*+s_*)}{\|X\|
^2\overline\psi(S_0)^2} \lesssim\frac{\log p (D_*+s_*)}{\overline\psi(S_0)^2}\\
&=:&R_*^2,
\end{eqnarray*}
this tends to zero. Thus we have proved that for some sufficiently
large constant~$M$,
\[
\E_{\be^0}\Pi \bigl(\be: \bigl\|X\bigl(\be-\be^0\bigr)
\bigr\|_2\geq 4\bigl\|X\bigl(\be^*-\be^0\bigr)\bigr\|_2+MR_*|
Y \bigr)\ra0.
\]
%

\emph{Proof of second assertion}.
Similar to \eqref{EqInvokeCompCond},
\begin{eqnarray*}
 &&\ola\bigl\|\be-\be^0\bigr\|_1 \\
 &&\qquad\leq\ola\bigl\|\be^*-\be^0
\bigr\|_1+ \frac{1}2 \bigl\| X\bigl(\be-\be^*\bigr)
\bigr\|_2^2 +\frac{\ola^2|S_{\be-\be^*}|}{2\|X\|^2\overline\psi(S_0)^2}
\\
&&\qquad\leq\bigl\|X\bigl(\be-\be^0\bigr)\bigr\|_2^2 +\ola\bigl\|
\be^*-\be^0\bigr\|_1+ \bigl\|X\bigl(\be^*-\be ^0\bigr)
\bigr\|_2^2 +\frac{\ola^2|S_{\be-\be^*}|}{2\|X\|^2\overline\psi(S_0)^2}.
\end{eqnarray*}
The claim follows now from the first assertion.



\emph{Proof of third assertion}. Note that $\|X(\be-\be^0)\|_2\geq
\widetilde\phi(|S_{\be-\be^0}|)\|X\|\|\be-\be^0\|_2\geq
\widetilde\psi(S_0)\|X\|\|\be-\be^0\|_2$. Now, the proof follows
from the first assertion.
\end{pf*}

\begin{pf*}{Proof of Theorem~\ref{thmm.BvM_type}}
The total variation distance between a probability measure $\Pi$
and its renormalized restriction $\Pi_A(\cdot):=\Pi(\cdot\cap
A)/\Pi(A)$ to a set
$A$ is bounded above by $2\Pi(A^c)$. We apply this to both the posterior
measure $\Pi(\cdot| Y)$ and the approximation $\Pi^\infty(\cdot| Y)$,
with the set
\[
A:= \biggl\{\be: \bigl\|\be-\be^0\bigr\|_1\leq \frac{Ms_0\sqrt{\log p}}{\|X\|
\overline\psi(S_0)^2\phi(S_0)^2}
\biggr\},
\]
where $M$ is a sufficiently large constant. By Theorem~\ref{TheoremRecovery}
the probability $\Pi(A| Y)$ tends to one under $\P_{\beta^0}$, and
at the end
of this proof we show that $\Pi^\infty(A| Y)$ tends to one as well.
Hence it suffices to prove
Theorem~\ref{thmm.BvM_type} with $\Pi(\cdot| Y)$ and
$\Pi^\infty(\cdot| Y)$ replaced by their renormalized restrictions
to $A$.

The measure $\Pi_A^\infty(\cdot| Y)$ is by its definition a
mixture over measures corresponding to models $S\in\cS_0$. By
Theorems~\ref{thmm.mod_sel} and~\ref{TheoremRecovery}
the measure $\Pi_A(\cdot| Y)$ is asymptotically concentrated on these models.
If $(\tilde v_S)$ is the renormalized restriction of a probability
vector $(v_S)$ to a set $\cS_0$,
then, for any probability measures $\Pi_S$,
\[
\biggl\|\sum_S \tilde v_S\Pi_S-
\sum_Sv_S\Pi_S
\biggr\|_{\TV} \le\bigl\|(\tilde v_S)-(v_S)
\bigr\|_{\TV}\le2\sum_{S\notin\cS_0}v_S,
\]
by the preceding paragraph. We infer that we can make a further reduction
by restricting and renormalizing the mixing weights of $\Pi(\cdot| Y)$
to $\cS_0$. More precisely, define probability measures by
\begin{eqnarray*}
\Pi^{(1)}(B| Y) &\propto&\sum_{S \in\mathcal{S}_0}
\frac{\pi_p(s)}{{p\choose s}} \biggl(\frac{\la}2 \biggr)^s \int
_{(B \cap A)_S} e^{-({1}/2)\|Y-X_S\be_S\|_2^2} e^{-\la\|\be
_S\|_1} \,d\be_S,
\\
\Pi^{(2)}(B| Y) &\propto&\sum_{S \in\mathcal{S}_0}
\frac{\pi_p(s)}{{p\choose s}} \biggl(\frac{\la}2 \biggr)^s \int
_{(B \cap A)_S} e^{-({1}/2)\|Y-X_S\be_S\|_2^2} e^{-\la\|\be
^0\|_1} \,d\be_S.
\end{eqnarray*}
Then it suffices to show that $\E_{\be^0}\|\Pi^{(1)}(\cdot| Y)-\Pi
^{(2)}(\cdot| Y)\|_{\TV}\ra0$.
(The factor $e^{-\la\|\be^0\|_1}$ in the second formula cancels in
the normalization,
but is inserted to connect to the remainder of the proof.)

For any sequences of measures $(\mu_S)$ and $(\nu_S)$, we have
\[
\biggl\|\frac{\sum_S\mu_S}{\|\sum_S\mu_S\|_{\TV}}
-\frac{\sum_S\nu_S}{\|\sum_S\nu_S\|_{\TV}} \biggr\|_{\TV}
\le\frac{2\sum_S\|\mu_S-\nu_S\|_{\TV}}{\|\sum_S\mu_S\|_{\TV}}
\le2\sup_S \biggl\|1-\frac{d\nu_S}{d\mu_S} \biggr\|_\infty
\]
if $\nu_S$ is absolutely continuous with respect to $\mu_S$ with
density $d\nu_S/d\mu_S$,
for every~$S$. It follows that
\begin{eqnarray*}
\bigl\|\Pi^{(1)}(\cdot| Y)-\Pi^{(2)}(\cdot| Y)\bigr\|_{\TV} &
\leq&2\max_{S\in\cS_0} \sup_{\be\in A}\bigl |
e^{\la\|\be_S\|
_1-\la\|\be^0\|_1}-1 \bigr|
\\
&\le&2\max_{S\in\cS_0} \sup_{\be\in A}
e^{\la\|\be_S-\be^0\|_1} \la\bigl\|\be_S-\be^0\bigr\|_1.
\end{eqnarray*}
This tends to zero by the definition of $A$ and the assumptions on $\be^0$.

Finally we show that $\Pi^\infty(A| Y)\ra1$. For $\Lambda_{n,\be
,\be^0}$, the
likelihood ratio given in~\eqref{eq.lr_representation}, we have
\begin{eqnarray}
\Pi^\infty\bigl(A^c| Y\bigr)= \frac{\int_{A^c} \Lambda_{n,\be,\be^0}(Y)
\,dU(\be)}{
\int\Lambda_{n,\be,\be^0}(Y) \,dU(\be)} \nonumber\\
\eqntext{\mbox{for }\displaystyle  dU(\be)= \sum_{S \in\cS_0} \frac{\pi_p(s)}{{p\choose s}} \biggl(
\frac{\la}2 \biggr)^s \,d\be_S\otimes
\delta_{S^c}.}
\end{eqnarray}
By \eqref{eq.lr_representation} the denominator in $\Pi^\infty(\cdot
| Y)$ satisfies
\begin{eqnarray*}
&& \int\Lambda_{n,\be,\be^0}(Y) \,dU(\be)
\\
&&\qquad\geq\frac{\pi_p(s_0)}{{p\choose s_0}}
\biggl(\frac{\la}2 \biggr)^{s_0} \int
e^{-({1}/2) \|X(\be_{S_0}-\be_{S_0}^0)\|_2^2+(Y-X\be^0)^t
X(\be_{S_0}-\be_{S_0}^0)} \,d\be_{S_0}
\\
&&\qquad\geq\frac{\pi_p(s_0)}{{p\choose s_0}}
 \biggl(\frac{\la}2 \biggr)^{s_0} \int
e^{-({1}/2 )\|Xb_{S_0}\|_2^2} \,db_{S_0} = \frac{\pi_p(s_0)}{{p\choose s_0}} \biggl(
\frac{\la}2 \biggr)^{s_0} \frac{(2\pi)^{s_0/2}}{
|\Gamma_{S_0}|^{1/2}},
\end{eqnarray*}
where $\Gamma_S=X_S^tX_S$, and
for the second inequality we use Jensen's inequality similarly  as in the proof
of Lemma~\ref{lem.expl_bd_for_denom}.

Using H\"older's inequality $|\al^t\be|\le\|\al\|_\infty\|\be\|_1$,
we see that on the event $\cT_{0}=\{\|X^t(Y-X\be^0)\|_\infty\le\ola
\}$,
\begin{eqnarray*}
\bigl(Y-X\be^0\bigr)^tX\bigl(\be-\be^0
\bigr)&\le&\ola\bigl\|\be-\be^0\bigr\|_1
\\
&\le&2\frac{\ola\|X(\be-\be^0)\|_2|S_{\be-\be^0}|^{1/2}}{\|X\|
\overline\phi(|S_{\be-\be^0}|)} -\ola\bigl\|\be-\be^0\bigr\|_1
\\
&\le&\frac{1}2 \bigl\|X\bigl(\be-\be^0\bigr)\bigr\|_2^2+
\frac{2\ola^2|S_{\be-\be^0}|}{
\|X\|^2\overline\phi(|S_{\be-\be^0}|)^2} -\ola\bigl\|\be-\be^0\bigr\|_1.
\end{eqnarray*}
Since $\ola(|S_{\be-\be^0}|)\geq\overline\psi(|S_0|)$ for every
$S_\be\in\cS_0$,
it follows that on $\cT_0$ the numerator in $\Pi^\infty(A^c| Y)$ is
bounded above by
\begin{eqnarray*}
&&e^{({2 \ola^2|S_{\be-\be^0}|}/{(\|X\|^2\overline\psi(S_0)^2)})
- ({\ola Ms_0 \sqrt{\log p}}/{(2\|X\|\overline\psi(S_0)^2\phi(S_0)^2)})}
\int e^{-({1}/2)\ola\|\be-\be^0\|_1} \,dU(\be)
\\
&& \qquad\leq e^{({8 |S_{\be-\be^0}|\log p}/{(\overline\psi(S_0)^2)}) -
({Ms_0 \log p}/{(2\overline\psi(S_0)^2\phi(S_0)^2)})} \sum_{s=0}^p
\pi_p(s)4^s.
\end{eqnarray*}
It follows that $\Pi^\infty(A^c| Y)$ is bounded above by
\[
\frac{{p\choose s_0}}{\pi_p(s_0)} \biggl(\frac{2}{\la} \biggr)^{s_0}
\frac
{|\Gamma_{S_0}|^{1/2}}{(2\pi)^{s_0/2}}
e^{({8 |S_{\be-\be^0}|\log p}/{\overline\psi(S_0)^2}) -
({Ms_0 \log p}/{(2\overline\psi(S_0)^2\phi(S_0)^2)})} \sum_{s=0}^p
\pi_p(s)4^s.
\]
By Jensen's inequality applied to the logarithm
$|\Gamma_S|\leq(s^{-1}\Tr(\Gamma_S))^s\leq\|X\|^{2s}$, and hence
$|\Gamma_S|^{1/2}/\la^{s}\le p^{s}$, by \eqref{eq.lambda_cond}.
The prior mass $\pi_p(s)$ can be bounded below by powers of $p^{-s}$ by
\eqref{assump.on_the_dim_prior}. This shows that the display tends to zero
for sufficiently large $M$.
\end{pf*}

\begin{pf*}{Proof of Theorem~\ref{TheoremSelectionNoSupersets}}
Let $\Sigma$ be the collection of all sets $S\in\cS_0$ such that
$S\supset S_0$ and
$S\neq S_0$. In view of Theorem~\ref{thmm.BvM_type} it suffices to
show that
$\Pi^\infty(\be: S_\be\in\Sigma| Y)\ra0$.

Note that due to $A_4>1$, any set in $S\in\cS_0$ has cardinality
smaller $6s_0$. By \eqref{EqDefWeightsw}, with $\Gamma_S=X_S^tX_S$,
\begin{eqnarray*}
 \Pi^\infty(\be: S_\be\in\Sigma| Y)& \leq&\sum
_{S\in\Sigma} \widehat w_S
\\
& \le&\sum_{s=s_0+1}^{6s_0}
{
\frac{\pi_p(s){p\choose s_0}{p-s_0\choose s-s_0} }{\pi
_p(s){p\choose s}}} \mathop{\max_{S\in\Sigma,}}_{{|S|=s}}
\frac{|\Gamma
_{S_0}|^{1/2}}{|\Gamma_{S}|^{1/2}} \biggl(\la\sqrt{\frac{\pi
}{2}} \biggr)^{s-s_0}
\\
&&{}\times e^{({1}/2 )\|X\widehat\be_{(S)}\|_2^2-({1}/2) \|X\widehat\be
_{(S_0)}\|_2^2}. 
\end{eqnarray*}
We shall show below that the factors on the right-hand side can be
bounded as follows:
for any fixed $r>2$,
%
\begin{eqnarray}
\label{EqSeparateInequalities} 
\la^{s-s_0} |
\Gamma_{S_0}|^{1/2}|\Gamma_{S}|^{-1/2} &
\leq&(4\sqrt{\log p})^{s-s_0} \widetilde\psi(S_0)
^{s_0-s},
\\
\quad\P \bigl(\|X_S\widehat\be_{(S)}\|_2^2-
\|X_{S_0}\widehat\be _{(S_0)}\|_2^2 &\le&
r (s-s_0)\log p, \mbox{ for all }S\in\Sigma \bigr)\ra1. \label{EqSeparateInequalitiesTwo}
\end{eqnarray}
Combining these estimates with assumption \eqref{assump.on_the_dim_prior}
shows that for $\cT$, the event in the second relation,
\[
\Pi^\infty(\be: S_\be\in\Sigma| Y)\mathbb{I}_{\cT}
\leq\sum_{s=s_0+1}^{6s_0} \bigl(A_1p^{-A_4}
\bigr)^{s-s_0} \pmatrix{s\cr s_0} \biggl(\frac{\sqrt{8\pi\log p}}{\widetilde\psi(S_0)}
\biggr)^{s-s_0} p^{r(s-s_0)/2}.
\]
For $s_0\le p^a$ we have ${s\choose s_0}={s\choose s-s_0}\le
s^{s-s_0}\le(6p^a)^{s-s_0}$.\vspace*{1.5pt}
Thus the expression tends to zero if $a-A_4+r/2<0$. Since $r$ can be chosen
arbitrarily close to $2$, this translates into $a<A_4-1$.

To prove bound \eqref{EqSeparateInequalities}, we apply the
interlacing theorem to
the principal submatrix $\Gamma_{S_0}$ of $\Gamma_S$ to see that
$\la_j(\Gamma_{S_0})\leq\la_j(\Gamma_S)$, for $j=1,\ldots,s_0$,
where $\la_1\ge\la_2\ge\cdots$
denote the eigenvalues in decreasing order, whence
\begin{eqnarray*}
|\Gamma_{S_0}|=\prod_{j=1}^{s_0}
\la_j(\Gamma_{S_0}) &\leq&\prod
_{j=1}^{s_0} \la_j(\Gamma_S)
\leq\la_{\min}(\Gamma _S)^{s_0-s} |
\Gamma_S|
\\
&\leq&\bigl(\widetilde\phi\bigl(|S|\bigr) \|X\|\bigr)^{2(s_0-s)} |\Gamma_S|.
\end{eqnarray*}
Assertion \eqref{EqSeparateInequalities} follows upon combining this
with \eqref{eq.lambda_cond}.

To bound the probability of the event $\cT$ in \eqref
{EqSeparateInequalitiesTwo},
we note that by the projection property of the least squares estimator,
for $S\supset S_0$ the difference $\|X_S\widehat\be_{(S)}\|_2^2-\|
X_{S_0}\widehat\be_{(S_0)}\|_2^2$ is the
square length of the projection of $Y$ onto the orthocomplement of the
range of $X_{S_0}$
within the range of $X_S$, a subspace of dimension $s-s_0$.
Because the mean $X\be^0$ of $Y=X\be^0+\ep$ is
inside the smaller of these ranges, it cancels under the projection,
and we may use the projection of the standard normal vector $\ep$ instead.
Thus the square length possesses a chi-square distribution with $s-s_0$
degrees of freedom.
There are $N={p-s_0 \choose s-s_0}$ models $S\in\Sigma$ that give rise
to such a chi-square distribution. Since $\log N\le(s-s_0) \log p\vee
1$, we can apply
Lemma~\ref{lem-chisq} with $d=s-s_0$ to give that
$\P(\cT^c)$ is bounded above by
$\sum_{s>s_0} {p-s_0\choose s-s_0}^{-(r-2)/4}e^{c(s-s_0)}$.
This tends to zero as $p\ra\infty$,
due to ${p-s_0\choose s-s_0}\ge(p-s)^{s-s_0}\ge(p/2)^{s-s_0}$, where
the last inequality follows from $s_0/p\leq s_0\la/\|X\|\rightarrow0$.
\end{pf*}

\begin{lem} \label{lem-chisq}
For every $r>2$, there exists a constant $c$ independent of $N\ge2$
and $d\ge1$
such that for any variables $Q_1,\ldots,Q_N$ that are marginally
$\chi^2(d)$ distributed,
\[
\mathbb{P} \Bigl( \max_{1\le i\le N} Q_i > r\log N
\Bigr) \le \biggl({\frac{1}N} \biggr)^{(r-2)/4}e^{cd}.
\]
\end{lem}

\begin{pf}
By Markov's inequality, for any $u>0$,
\[
\P \Bigl( \max_{1\le i\le N} Q_i > r\log N \Bigr)  \le
e^{-ur\log N} \E\max_{1\le i\le N} e^{uQ_i} \le
N^{-ur} N \sqrt{1-2u}^{-d}.
\]
The results follows upon choosing $u=1/4+1/(2r)$, giving $ur-1=(r-2)/4$
and $1-2u=1/2-1/r$.
\end{pf}

\begin{pf*}{Proof of Theorem~\ref{TheoremSelection}}
\emph{Proof of first two assertions}.
Because $\|\be_{S_0}-\be^0\|_1\le\|\be-\be^0\|_1$, the posterior
probability of the
set
\[
\biggl\{\be: \bigl\|\be_{S_0}-\be^0\bigr\|_1>
\frac{M}{\overline\psi(S_{0})^2} \frac{|S_{0}| \sqrt{\log p}}{\|X\|\phi(S_{0})^2} \biggr\}
\]
tends to zero by Theorem~\ref{thmm.pred_and_l1}. This implies the
first assertion.
The second assertion follows similarly from the second assertion
of Theorem~\ref{thmm.pred_and_l1}.

\emph{Proof of third assertion}.
First we prove that the largest coefficient in absolute value,
say $\be^0_{m}$, is selected by the posterior if this is above the threshold.
By Theorem~\ref{thmm.BvM_type} it is enough to show
that $\E_{\be^0}\Pi^\infty (\be: m\in S_\be| Y )\ra1$.
For any given set $S$ with $m\notin S$, let $S_m:=S \cup\{ m \}$ and $s=|S|$.
Then
\[
\Pi^\infty(\be: m\notin S_\be| Y)  = \sum
_{S\in\cS_0: m\notin S} \widehat w_S.
\]
We shall bound this further by showing that $\widehat w_S\ll\widehat
w_{S_m}$, for every $S$
in the sum. The quotient of these weights is equal to
\begin{eqnarray*}
\frac{\widehat w_{S_m}}{\widehat w_S} & = &\la\sqrt{\frac{\pi}{2}} \frac{\pi_p(s+1)}{\pi_p(s)}
\frac
{{p\choose s}}{{p\choose s+1}}
 \frac{ |\Gamma_{S}|^{1/2}}{|\Gamma_{S_m}|^{1/2}}
 e^{({1}/2) \|X_{S_m}\widehat\be_{(S_m)}\|_2^2-({1}/2) \|X_S\widehat
\be_{(S)}\|_2^2}
\\
& \geqa&\la p^{-A_3} \frac{s+1}{p-s} \frac{ |\Gamma
_{S}|^{1/2}}{|\Gamma_{S_m}|^{1/2}}
e^{({1}/2) \|X_{S_m}\widehat\be_{(S_m)}\|_2^2-({1}/2) \|X_S\widehat
\be_{(S)}\|_2^2},
\end{eqnarray*}
in view of \eqref{assump.on_the_dim_prior}.
By the interlacing theorem, the eigenvalues $\overline\la_i$
in increasing order of the matrices $\Gamma_S$ and
$\Gamma_{S_m}$ satisfy $\overline\la_i(\Gamma_{S_m})\le
\overline\la_i(\Gamma_{S})\le\overline\la_{i+1}(\Gamma_{S_m})$,
for any $1\le i\le s$. This
implies that $|\Gamma_{S}|/|\Gamma_{S_m}|\ge\overline\la
_{s+1}(\Gamma_{S_m})^{-1}$. Since
$\|X\be\|_2 \le\|X\|\|\be\|_1\le\sqrt{|S_\be|}\|X\|\|\be\|_2$,
for any $\be$, the largest
eigenvalue $\overline\la_{s+1}(\Ga_{S_m})$ is at most $(s+1)\|X\|
^2$. Combining
this with \eqref{eq.lambda_cond},
we conclude that the preceding display is bounded below by
\begin{eqnarray*}
&&\frac{\la}{\|X\|} p^{-A_3-1} e^{({1}/2) \|X_{S_m}\widehat\be
_{(S_m)}\|_2^2-({1}/2) \|X_S\widehat\be_{(S)}\|_2^2} \\
&&\qquad\ge
p^{-A_3-2} e^{({1}/2) \|X_{S_m}\widehat\be_{(S_m)}\|_2^2-({1}/2) \|X_S\widehat\be_{(S)}\|_2^2}.
\end{eqnarray*}
By definition of the least squares estimator, the difference
of the square norms in the exponent is the square length of the
projection of $Y=X\be^0+\ep$ onto the orthocomplement $F_S$ of the
range of $X_S$ in
the range of $X_{S_m}$, the one-dimensional space spanned by the vector
$X_{m}-P_SX_{m}$, where $P_S$ denotes the projection onto the range of
$X_S$. If,
with an abuse of notation, $P_{F_S}$ is the projection onto $F_S$, then
%
\begin{eqnarray}\label{EqSizeProjection}
\qquad\|X_{S_m}\widehat\be_{(S_m)}\|_2^2
- \|X_S\widehat\be_{(S)}\|_2^2 & =&
\| P_{F_S} Y \|_2^2 \ge\frac{1}2
\bigl\|P_{F_S} X\be^0\bigr\|_2^2-\|
P_{F_S}\ep\|_2^2
\nonumber
\\[-8pt]
\\[-8pt]
\nonumber
&=&\frac{\langle X\be^0,X_m-P_SX_m\rangle^2}{2\|X_m-P_SX_m\|_2^2}
 -\frac{\langle\ep,X_m-P_SX_m\rangle^2}{\|X_m-P_SX_m\|_2^2}.
\end{eqnarray}
We shall show that the first term on the right
is large if $|\be_m^0|$ is large, and the second is small with large
probability.

We start by noting that for $j\notin S$ and any $S$,
%
\begin{eqnarray}
\label{EqProjNormEstimate} \|P_S X_{ j}\|_2^2&=&
\bigl\langle X_{j}, X_S\Gamma_S^{-1}X_S^tX_{j}
\bigr\rangle \le\frac{1}{\widetilde\phi(s)^2\|X\|^2}\bigl\|X_S^tX_j
\bigr\|_2^2
\nonumber
\\[-8pt]
\\[-8pt]
\nonumber
&=&\frac{1}{\widetilde\phi(s)^2\|X\|^2}\sum_{i\in S}\bigl(X^tX
\bigr)_{i,j}^2 
\le\frac{s\mc(X)^2\|X\|^2}{\widetilde\phi(s)^2}.
\end{eqnarray}
%
It follows from the definitions that
$\widetilde\phi(1)\|X\|\le\|X_j\|\le\|X\|$, for every $j$.
Combined, this shows that $\|X_j-P_SX_j\|_2\ge3\|X\|\widetilde\phi
(1)/4$ if
$\sqrt s\mc(X)\le\widetilde\phi(s)\widetilde\phi(1)/4$.

We write $X\be^0=X_m\be_m^0+X_{-m}\be_{-m}^0$, for $X_{-m}$ the
matrix obtained
by removing the column $X_m$ from $X$, and split the first inner product
in \eqref{EqSizeProjection} in the two parts
\begin{eqnarray*}
\bigl|\bigl\langle X_m\be_m^0,
X_m-P_SX_m\bigr\rangle\bigr| &=&\bigl|
\be_m^0\bigr| \|X_m-P_SX_m
\|_2^2,
\\
\bigl|\bigl\langle X_{-m}\be_{-m}^0,
X_m-P_SX_m\bigr\rangle\bigr| &=&\biggl|\sum
_{j\neq m}\be_j^0 \langle
X_{j}-P_SX_j, X_m-P_SX_m
\rangle\biggr|
\\
&\le&\sum_{j\neq m, j\notin S}\bigl|\be_j^0\bigr|
\bigl(\mc(X)\|X\|^2+\| P_SX_j
\|_2 \|P_SX_m\|_2 \bigr)
\\
&\le& s_0\bigl |\be_m^0\bigr| \biggl(\mc(X)\|X
\|^2+\frac{s\mc(X)^2\|X\|
^2}{\widetilde\phi(s)^2} \biggr),
\end{eqnarray*}
using that $X_j-P_SX_j=0$ if $j\in S$,
the definition of $\mc(X)$ to bound $\langle X_j,X_m\rangle$, the
Cauchy--Schwarz inequality on $\langle P_SX_j, X_m\rangle=\langle
P_SX_jz P_SX_m\rangle$ and \eqref{EqProjNormEstimate}. Putting the
estimates together we find that
for $(s_0\vee s)\mc(X)\le\widetilde\phi(s)\widetilde\phi(1)/4$,
\[
\bigl\|P_{F_S}X\be^0\bigr\|_2 \ge\bigl|\be_m^0\bigr|
\|X\| \widetilde\phi(1)\tfrac{1}4. 
\]
We can split the random inner product in \eqref{EqSizeProjection}
in the two parts $\langle\ep,X_m\rangle$ and $\langle\ep
,P_SX_m\rangle$.
For $\sqrt s\mc(X)\le\widetilde\phi(s)\widetilde\phi(1)/2$,
\[
\|P_{F_S}\ep\|_2 \le\frac{|\langle\ep,X_m\rangle|}{3\|X\|\widetilde\phi(1)/4} +
\frac{|\langle\ep,P_SX_m\rangle|}{3\|X\|\widetilde\phi(1)/4}.
\]
Each variable $\langle\ep, v\rangle$ is normally distributed with
mean zero and
variance $\|v\|_2^2$, for any $v\in\RR^n$. When $m$ varies over
$1,\ldots, p$
and $S$ over all subsets of size $s$ that do not contain $m$, there are
$p$ possible
variables in the first term and $p{p-1 \choose s}$ possible variables
in the second.
For $\widetilde\phi(s)\ge\widetilde\psi(S_0)\ge c_0$ the variances
of the variables in the two terms are of the orders
$1/c_0^2$ and $s\mc(X)^2/c_0^4$, respectively.
Therefore the means of the two suprema are of the orders $\sqrt{\log
p}$ and
$\sqrt{\log{p \choose s}}s^{1/2}\mc(X)\le\sqrt{\log p}$, respectively,
if $s\mc(X)\le1$.
With probability $O(p^{-\mu})$ these variables do not exceed a multiple
of their means.

We conclude that for $(s_0\vee s)\mc(X)\le\widetilde\phi
(s)\widetilde\phi(1)/4$ and
$\widetilde\phi(s)\ge c_0$, the left-hand side of \eqref{EqSizeProjection}
is, with probability tending to one, bounded below by $\|X\|^2(\be
_m^0)^2c_0^2/16-O(\log
p)$,
whence for $|\be^0_m|\ge M\sqrt{\log p}/\|X\|$ for large $M$,
uniformly in $S,m$,
\[
\frac{\widehat w_{S_m}}{\widehat w_S}\ge p^{-A_3-2}e^{c M^2\log p}\ge p^\mu,
\]
for $\mu>0$ as large as desired (depending on $M$) and $c$ a suitable
positive constant.
So, with overwhelming probability, 
\[
\Pi^\infty(\be: m\notin S_\be| Y) \le p^{-\mu}
\sum_{S\in\cS_0: m\in S} \widehat w_S\le
p^{-\mu}.
\]
Thus $\E_{\be^0}\Pi^\infty (m\notin S| Y )\ra0$ at the
order $p^{-\mu}$.

Next, for $\be_{m_2}$ the second largest coefficient, we consider $\Pi
^\infty (m_2\notin
S| m_1\in S, Y )$. By reasoning similar to the preceding, we
show that the index $m_2$ is included asymptotically, etc.
\end{pf*}

\section*{Acknowledgments}
We thank an Associate Editor and four referees for valuable comments.
We are also grateful to Amandine Schreck for helpful discussions.

\begin{supplement}[id=suppA]
\stitle{Bayesian linear regression with sparse priors\\}
\slink[doi]{10.1214/15-AOS1334SUPP} 
\sdatatype{.pdf}
\sfilename{aos1334\_supp.pdf}
\sdescription{In the supplement
we state a Bernstein--von Mises type result for large lambda and give
the remaining proofs.}
\end{supplement}






\printaddresses
\end{document}